\documentclass[12pt]{amsart}

\usepackage{amsmath,amssymb,amsthm}
\usepackage{epsf}
\usepackage{epsfig}
\usepackage[T2A]{fontenc}
\usepackage[cp1251]{inputenc}
\usepackage[english]{babel}

\newcommand{\RR}{\mathbb R}
\newcommand{\CC}{\mathbb C}

\renewcommand{\Re}{\mathop{\rm Re}\nolimits}
\renewcommand{\Im}{\mathop{\rm Im}\nolimits}
\newcommand{\trace}{\mathop{\mathrm{trace}}\nolimits}

\newcommand{\ii}{\mathrm{i}}

\newcommand{\manifold}[1]{\mathcal{#1}}
\newcommand{\M}{\manifold{M}}
\newcommand{\D}{\manifold{D}}

\newcommand{\vect}[1]{\mathrm{#1}} %{\mathbf{#1}}
\newcommand{\x}{\vect{x}}
\newcommand{\y}{\vect{y}}
\newcommand{\va}{\vect{a}}
\newcommand{\vb}{\vect{b}}
\newcommand{\vc}{\vect{c}}
\newcommand{\vd}{\vect{d}}
\newcommand{\vX}{\vect{X}}
\newcommand{\vY}{\vect{Y}}
\newcommand{\vH}{\vect{H}}
\newcommand{\n}{\vect{n}}

\newtheorem{thm}{Theorem}[section]

\theoremstyle{definition}
\newtheorem{defn}[thm]{Definition}

\theoremstyle{remark}
\newtheorem{rem}[thm]{Remark}

\newcommand{\prfend}{\hfill $\blacksquare$ \bigskip}

\newcommand{\ds}{\displaystyle}

\makeatletter
\@addtoreset{equation}{section}
\makeatother

\begin{document}

\title [Canonical Weierstrass representations for minimal space-like surfaces]
{Canonical Weierstrass representations for minimal space-like surfaces in $\RR^4_1$}

\author{Georgi Ganchev and Krasimir Kanchev}

\address{Bulgarian Academy of Sciences, Institute of Mathematics and Informatics,
Acad. G. Bonchev Str. bl. 8, 1113 Sofia, Bulgaria}
\email{ganchev@math.bas.bg}%

\address {Department of Mathematics and Informatics, Todor Kableshkov University of Transport,
158 Geo Milev Str., 1574 Sofia, Bulgaria}%
\email{kbkanchev@yahoo.com}%


\begin{abstract}
A space-like surface in Minkowski space-time is minimal if its mean curvature vector field
is zero. Any minimal space-like surface of general type admits special isothermal parameters -
canonical parameters. For any minimal surface of general type parameterized by canonical
parameters we obtain Weierstrass representations - canonical Weierstrass representations
via two holomorphic functions. We find the expressions of the Gauss curvature and the normal
curvature of the surface with respect to this pair of holomorphic functions. We find the
relation between two pairs of holomorphic functions generating one and the same minimal
space-like surface of general type. The canonical Weierstrass formulas allow us to establish geometric
correspondence between minimal space-like surfaces of general type and classes of pairs
of holomorphic functions in the Gauss plane.
\end{abstract}

\thispagestyle{empty}

\maketitle

\section{Introduction}

A two-dimensional surface $\M$ in the four-dimensional Minkowski space-time $\RR^4_1$ is said
to be \emph{space-like} if the induced metric on the tangential space at any point of $\M$ is
positive definite. If $\M$ is a space-like surface in $\RR^4_1$, we denote by $T_p(\M)$ and
$N_p(\M)$ the tangential space and the normal space at a point $p \in \M$, respectively.
The flat Levi-Civita connection on $\RR^4_1$ is denoted by $\nabla$. Then the second fundamental
tensor $\sigma$ of $\M$ is given by
$$\sigma (X,Y) = (\nabla_X Y)^\bot ;\quad X,Y \; \text{tangent vectors to}\; \mathcal{M}\;
\text{at a point} \; p\in \M.$$

The space-like surface $\M$ is minimal if its mean curvature vector field $H=\frac{1}{2}\trace \sigma$ is zero,
i.e. $H =0$.

A general approach to Weierstrass representations of minimal space-like surfaces  in $\RR^4_1$
was given in \cite{dM} and \cite{E-1}.

In \cite{A-P-1} minimal space-like surfaces in $\RR^4_1$ were studied with respect to special isothermal
parameters and a fundamental theorem of Bonnet type in terms of the Gauss curvature $K$
and the normal curvature $\varkappa$ was proved.

The question when a complete minimal space-like surface is a plane was studied in \cite{E-R-1}.

In this paper we consider canonical Weierstrass representations for minimal space-like surfaces in $\RR^4_1$.

A point $p\in\M$ is said to be \emph{degenerate}, if the set $\{\sigma(\vX,\vY);\ \vX \in T_p(\M),\vY \in T_p(\M) \}$,
is contained in one of the two light-like one-dimensional subspaces of $N_p(\M)$.

We call a minimal space-like surface, free of degenerate points, a {\it minimal space-like surface of general type}.

Let $(\mathcal{M},\text{x}(u, v))$ be a space-like surface in $\RR^4_1$,
parameterized by isothermal coordinates $(u, v)$. In isothermal coordinates the space-like surface $\M$
is minimal if and only if the position vector function $\text{x}(u,v)$ is harmonic.

We describe the properties of minimal surfaces in terms of the complex vector function $\Phi(t)=\x_u-\ii\x_v , \;
t=u+\ii v $.

The standard Weierstrass representations for minimal space-like surfaces are in terms of three holomorphic functions.

Using special isothermal parameters (canonical parameters) on a minimal space-like surface
of general type, we obtain canonical Weierstrass representations in terms of two holomorphic functions.

We call the isothermal parameters {\it canonical of the first type $($the second type$)$} if
${\Phi^\prime}^2=+1 \; ({\Phi^\prime}^2=-1)$. The special parameters, used in \cite{A-P-1} occur to be
canonical of the first type.

In Theorem \ref{can_coord_exist} we prove that:

{\it Any minimal space-like surface in $\RR^4_1$, free of degenerate points, admits locally canonical
coordinates of both types.}

In Theorem \ref{thm_Wcang} we prove the following statement.

{\it Any minimal space-like surface ${\M}$ of general type, parameterized by canonical coordinates
of the first type, has the following Weierstrass representation:
\begin{equation}\label{Wcang0}
\Phi: \quad
\begin{array}{rll}
\phi_1 &=& \ds\frac{\ii}{2}\; \ds\frac {g_1 g_2+1}{\sqrt{g'_1 g'_2}}\,,\\[6mm]
\phi_2 &=& \ds\frac{1}{2}\; \ds\frac {g_1 g_2-1}{\sqrt{g'_1 g'_2}}\,,\\[6mm]
\phi_3 &=& \ds\frac{1}{2}\; \ds\frac {g_1 + g_2}{\sqrt{g'_1 g'_2}}\,,\\[6mm]
\phi_4 &=& \ds\frac{1}{2}\; \ds\frac {g_1 - g_2}{\sqrt{g'_1 g'_2}}\,,\\
\end{array}
\end{equation}
where $(g_1,g_2)$ is a pair of holomorphic functions satisfying the conditions:
\begin{equation}\label{Wcang_cond0}
g'_1 g'_2 \neq 0; \quad g_1 \bar g_2 \neq -1.
\end{equation}

Conversely, if $(g_1,g_2)$ is a pair of holomorphic functions satisfying the conditions \eqref{Wcang_cond0},
then the formulas \eqref{Wcang0} generate a minimal space-like surface of general
type, parameterized by canonical coordinates of the first type.}

We call the representation of $\Phi$ in Theorem \ref{thm_Wcang} \emph{canonical Weierstrass representation}.

In terms of the above canonical representation the coefficients of the first fundamental
form are given by
$$E = G =\ds\frac{|1+g_1\bar g_2|^2}{4|g'_1 g'_2|}.$$

The Gauss curvature $K$ and the curvature of the normal connection $\varkappa$ (the normal curvature)
are given by
$$K = \Re \frac{-16|g'_1 g'_2|\, g'_1 \bar {g'_2}}
                          {|1 + g_1 \bar g_2|^2\; (1 + g_1 \bar g_2)^2}\,,\qquad
\varkappa = \Im \frac{-16|g'_1 g'_2|\, g'_1 \bar {g'_2}}
                          {|1 + g_1 \bar g_2|^2\; (1 + g_1 \bar g_2)^2}\;.$$

Theorem \ref{thm_Wcang} gives a representation of any minimal space-like surface of general type in terms of
two holomorphic functions. The following question arises naturally:

If $(g_1, g_2)$ and $(\hat g_1, \hat g_2)$ are two pairs of holomorphic functions generating
one and the same minimal space-like surface of general type, what is the relation between them?

We answer to this question in Theorem \ref{A-B_Cang}.

{\it Let $({\hat\M},\hat\x)$ and $({\M},\x)$ be two minimal space-like surfaces of general type, given by
the canonical Weierstrass representation of the type \eqref{Wcang0}. The following conditions are equivalent:
\begin{enumerate}
    \item $({\hat\M},\hat\x)$ and $({\M},\x)$ are related by a transformation in $\RR^4_1$ of the type:\\
    $\hat\x(t)=A\x(t)+\vb$, where $A \in \mathbf{SO}(3,1,\RR)$ and $\vb \in \RR^4_1$.
    \item The functions in the Weierstrass representations of $({\hat\M},\hat\x)$ and $({\M},\x)$ are related by the
    following equalities:
$$
\hat g_1 = \ds\frac{ag_1+b}{cg_1+d}\,; \quad \hat g_2 = \ds\frac{\phantom{-}\bar d g_2 - \bar c}{-\bar b g_2 + \bar a}\;,
\qquad a,b,c,d \in \CC, \; ad-bc\neq 0.
$$
\end{enumerate}}

\section{Preliminaries}
Let $\RR^4_1$ denote the standard Minkowski space-time. This is a
four-dimensional space endowed  with the indefinite dot product:
\begin{equation}\label{Mink}
\va\cdot \vb=a_1b_1+a_2b_2+a_3b_3-a_4b_4\; .
\end{equation}

If $\M$ is a two-dimensional manifold and $\x : \M \to \RR^4_1$ is an  immersion
of $\M$ into $\RR^4_1$, then we say that $\M$ is a (regular)
surface in $\RR^4_1$. We denote by $T_p(\M)$  the tangential space
to $\M$ at a point $p$ identifying $T_p(\M)$ with the
corresponding plane in $\RR^4_1$. $N_p(\M)$ will stand for the
normal space to $\M$ at the point $p$, which is the orthogonal
complement to $T_p(\M)$ in $\RR^4_1$. If the induced metric onto
$T_p(\M)$ is positive definite, the surface $\M$ is said to be
space-like. Then the induced metric onto the normal space
$N_p(\M)$ is of signature (1,1). The surface $\M$ with the induced
metric becomes a two-dimensional Riemannian space.

Let $E=\x_u^2$, $F=\x_u\cdot\x_v$ and $G=\x_v^2$ be the
coefficients of the first fundamental form on $\M$. The surface
$\M$ admits locally around any point $p \in \M$ isothermal
coordinates (parameters) $(u,v)$, which means that $E=G$ and
$F=0$. Together with the real coordinates $(u,v) \in \D$ we also
consider the complex coordinate $t=u+\ii v$, identifying $\RR^2$
with the complex plane $\CC$. Thus all functions defined around
$p$ can be considered as functions of the complex variable $t$.
Throughout this paper we consider isothermal coordinates $(u,v)$ on
$\M$.

We also consider the complexified tangential space $T_{p,C}(\M)$
and the complexified normal space $N_{p,C}(\M)$ at a point $p$ in
$\M$ as the corresponding 2-planes in $\CC^4$.

If $\va$ and $\vb$ are two vectors in $\CC^4$, then by $\va\cdot
\vb$ (or $\va\vb$) we denote the bilinear product in $\CC^4$,
which is the natural extension of the product in $\RR^4_1$ given
by \eqref{Mink}. Together with the bilinear product in $\CC^4$ we
also consider the indefinite Hermitian product of $\va$ and $\vb$,
given by
\[\va\cdot\bar \vb = a_1\bar b_1+a_2\bar b_2+a_3\bar b_3-a_4\bar b_4\; .\]
The square of $\va$ with respect to the bilinear product is
\[\va^2=\va\cdot \va = a_1^2+a_2^2+a_3^2-a_4^2\;\]
and the norm of $\va$ with respect to the Hermitian product is
\[\|\va\|^2=\va\cdot\bar \va = |a_1|^2+|a_2|^2+|a_3|^2-|a_4|^2\; .\]
The spaces $T_{p,C}(\M)$ and $N_{p,C}(\M)$ are closed with respect
to the complex conjugation and are orthogonal with respect to
both: bilinear and Hermitian product. Therefore we have the
following orthogonal decomposition
\[\CC^4 = T_{p,C}(\M) \oplus N_{p,C}(\M)\; .\]
For a given vector $\va \in \CC^4$ we denote by $\va^\top$ and
$\va^\bot$ the orthogonal projections of $\va$ into $T_{p,C}(\M)$
and $N_{p,C}(\M)$, respectively, i.e.
\[\va=\va^\top + \va^\bot\; .\]
The above decomposition does not depend on the bilinear or the
Hermitian dot product in $\CC^4$.

The second  fundamental form on $\M$ is denoted by $\sigma$. By
definition we have:
\[\sigma (\vX,\vY) = (\nabla_\vX \vY)^\bot \; , \]
where $\vX, \vY \in T(\M)$, and $\nabla$ is the canonical flat
connection in $\RR^4_1$.

Let $\vX_1$ and $\vX_2$ denote the unit tangent vectors to $\M$ at
a point $p$ having the same directions as  the coordinate vectors
$\x_u$ and $\x_v$, respectively, i.e.
\begin{equation*}%\label{\vX1\vX2}
\vX_1=\frac{\x_u}{\|\x_u\|}=\frac{\x_u}{\sqrt E};\quad
\vX_2=\frac{\x_v}{\|\x_v\|}=\frac{\x_v}{\sqrt G}=\frac{\x_v}{\sqrt
E}\; .
\end{equation*}

The mean curvature $\vH$ of $\M$ is the vector function
\[\vH = \frac{1}{2}\trace\sigma = \frac{1}{2}(\sigma(\vX_1,\vX_1)+\sigma(\vX_2,\vX_2)) \; . \]

A space-like surface $\M$ in $\RR^4_1$ is said to be minimal if
$\vH=0$ at any point of $\M$.

\section{The complex function $\Phi (t)$. }

Let $\M$ be a space-like surface in $\RR^4_1$, parameterized by
isothermal coordinates. The complex-valued vector function
$\Phi(t)$ on $\M$ with values in $\CC^4$ is defined by
\begin{equation}\label{Phi_def}
\Phi(t)=2\frac{\partial\x}{\partial t}=\x_u-\ii\x_v \; .
\end{equation}

The defining equality \eqref{Phi_def} implies that:
 \[ \Phi^2=(\x_u - \ii\x_v)^2=\x_u^2-\x_v^2-2\x_u \x_v \ii\; .\]
Then the following equalities are equivalent:
$$\Phi^2=0 \  \Leftrightarrow\  \begin{array}{l} \x_u^2-\x_v^2=0\\  \x_u \x_v=0 \end{array} \
\Leftrightarrow \  \begin{array}{l} E=\x_u^2=\x_v^2=G\\ F=0
\end{array}.$$ Hence, the parameters $(u,v)$ are isothermal if and
only if:
\begin{equation}\label{Phi2}
\Phi^2=0\; .
\end{equation}

For the norm of $\Phi$ we find
\[
\|\Phi\|^2=\Phi\bar{\Phi}=\x_u^2+\x_v^2=2E=2G.
\]
Therefore
\begin{equation}\label{EG}
E=G=\frac{1}{2}\|\Phi\|^2,\quad F=0
\end{equation}
and
\begin{equation}\label{Idt}
\mathbf{I}=\frac{1}{2}\|\Phi\|^2 (du^2 +
dv^2)=\frac{1}{2}\|\Phi\|^2|dt|^2 \; .
\end{equation}
From the above it follows that $\Phi$ satisfies the condition:
\begin{equation}\label{modPhi2}
\|\Phi\|^2> 0\; .
\end{equation}

Differentiating equality \eqref{Phi_def} and using that
$\frac{\partial}{\partial\bar t} \frac{\partial}{\partial t} =
\frac{1}{4}\Delta$, we find
\begin{equation}\label{dPhi_dbt}
\frac{\partial\Phi}{\partial\bar t}= \frac{\partial}{\partial\bar
t}\, \left(2\, \frac{\partial \x}{\partial t} \right)=
\frac{1}{2}\Delta \x \; ,
\end{equation}
where $\Delta$ denotes the Laplace operator.

The last formula implies that $\ds\frac{\partial\Phi}{\partial\bar
t}$ is a real vector function, i.e.
\begin{equation}\label{dPhi_dbt=dbPhi_dt}
\frac{\partial\Phi}{\partial\bar
t}=\frac{\partial\bar\Phi}{\partial t} \; .
\end{equation}

Thus, any space-like surface in $\RR^4_1$ parameterized by
isothermal coordinates, determines a function $\Phi$ given by
\eqref{Phi_def}, which satisfies the conditions:
\begin{equation}\label{Phi_cond}
\Phi^2=0, \quad \|\Phi\|^2> 0, \quad
\frac{\partial\Phi}{\partial\bar
t}=\frac{\partial\bar\Phi}{\partial t} \; .
\end{equation}
Conversely, any function $\Phi$ satisfying these three conditions
determines locally a space-like surface in isothermal coordinates up to a
translation.

The last assertion follows immediately from the fact that
\eqref{dPhi_dbt=dbPhi_dt} is the integrability condition for the
system
\begin{equation}\label{xuxv1}
\begin{array}{ll}
\x_u=\ \ \,\Re (\Phi),\\[2mm]
\x_v=-\Im (\Phi).
\end{array}
\end{equation}

Next we express the vectors $\x_u$, $\x_v$ and the second
fundamental form $\sigma$ of $\M$ by means of $\Phi$.

Taking into account \eqref{Phi_def} we have:
\begin{equation}\label{xuxv}
\begin{array}{ll}
\x_u=\ \ \,\Re (\Phi)=\ds\frac{1}{2}(\Phi+\bar\Phi),\\[3mm]
\x_v=-\Im (\Phi)=\ds\frac{\ii}{2}(\Phi-\bar\Phi).
\end{array}
\end{equation}

Equality \eqref{dPhi_dbt} implies that:
\begin{equation*}%\label{dPhi_dbt}
\left(\frac{\partial\Phi}{\partial\bar t}\right)^\bot=
\left(\frac{1}{2}\Delta \x \right)^\bot= \frac{1}{2}(\x_{uu}^\bot
+ \x_{vv}^\bot)= \frac{1}{2}({\nabla_{\x_u}^\bot \x_u} +
{\nabla_{\x_v}^\bot \x_v })= \frac{1}{2}(\sigma (\x_u,\x_u)+\sigma
(\x_v,\x_v)) .
\end{equation*}

Differentiating \eqref{Phi_def} with respect to $t$ we find
\begin{equation}\label{dPhi_dt}
\frac{\partial\Phi}{\partial t}= \frac{1}{2}(\x_{uu} - \x_{vv}) -
\ii \x_{uv} .
\end{equation}
and
\begin{equation}\label{dPhi_dt_bot}
\left(\frac{\partial\Phi}{\partial t}\right)^\bot=
\frac{1}{2}(\sigma (\x_u,\x_u)-\sigma (\x_v,\x_v)) - \ii \sigma
(\x_u,\x_v) .
\end{equation}
Therefore
\begin{equation}
\begin{array}{l}\label{sigma_uu_vv_uv}
\ds\sigma (\x_u,\x_u) = \Re\left(\frac{\partial\Phi}{\partial\bar
t}\right)^\bot +
                        \Re\left(\frac{\partial\Phi}{\partial     t}\right)^\bot;\\
\ds\sigma (\x_v,\x_v) = \Re\left(\frac{\partial\Phi}{\partial\bar
t}\right)^\bot -
                        \Re\left(\frac{\partial\Phi}{\partial     t}\right)^\bot;\\
\ds\sigma (\x_u,\x_v) = -\Im\left(\frac{\partial\Phi}{\partial
t}\right)^\bot.
\end{array}
\end{equation}

Finally we give transformation formulas for the function $\Phi$
under a change of the isothermal coordinates and under a motion in
$\RR^4_1$.

Let us consider the change of the isothermal coordinates given by
$t=t(s)$. Since the transformation of the isothermal coordinates
is conformal in $\CC$, then the function $t(s)$ is either
holomorphic or antiholomorphic. Denote by $\tilde\Phi(s)$ the
corresponding function in the new coordinates.

First, let us consider the holomorphic case. Taking into account
\eqref{Phi_def} we have:
\begin{equation*}
\tilde\Phi(s)=2\frac{\partial\x}{\partial s}=
2\frac{\partial\x}{\partial t}\frac{\partial t}{\partial s}+
2\frac{\partial\x}{\partial \bar t}\frac{\partial \bar t}{\partial
s}= 2\frac{\partial\x}{\partial t}\frac{\partial t}{\partial s} \;
.
\end{equation*}
Therefore, if the change $t=t(s)$ is holomorphic, then
\begin{equation}\label{Phi_s-hol}
\tilde\Phi(s)=\Phi(t(s)) \frac{\partial t}{\partial s} \; .
\end{equation}
In the antiholomorphic case we have similarly
\begin{equation*}
\tilde\Phi(s)=2\frac{\partial\x}{\partial s}=
2\frac{\partial\x}{\partial t}\frac{\partial t}{\partial s}+
2\frac{\partial\x}{\partial \bar t}\frac{\partial \bar t}{\partial
s}= 2\frac{\partial\x}{\partial \bar t}\frac{\partial \bar
t}{\partial s},
\end{equation*}
i.e.
\begin{equation}\label{Phi_s-antihol}
\tilde\Phi(s)=\bar\Phi(t(s)) \frac{\partial \bar t}{\partial s} \;
.
\end{equation}

Now let $(\M,\x)$ and $(\hat\M,\hat\x)$ be two surfaces in
$\RR^4_1$, parameterized by isothermal coordinates $t=u+\ii v$ in
one and the same domain $\D \subset \CC$. Suppose that
$(\hat\M,\hat\x)$ is obtained by $(\M,\x)$ by means of a motion in
$\RR^4_1$:
\begin{equation}\label{hat_M-M-mov}
\hat\x(t)=A\x(t)+\vb; \qquad A \in \mathbf{O}(3,1,\RR), \ \vb \in
\RR^4_1 \; .
\end{equation}
Differentiating \eqref{hat_M-M-mov} we find the relation between
the corresponding functions $\Phi$ and $\hat\Phi$:
\begin{equation}\label{hat_Phi-Phi-mov}
\hat\Phi(t)=A\Phi(t); \qquad A \in \mathbf{O}(3,1,\RR) \; .
\end{equation}
Conversely, if $\Phi$ and $\hat\Phi$ are connected by
\eqref{hat_Phi-Phi-mov}, then we have $\hat\x_u=A\x_u$ and
$\hat\x_v=A\x_v$, which implies \eqref{hat_M-M-mov}. Hence, the
relations \eqref{hat_M-M-mov} and \eqref{hat_Phi-Phi-mov} are
equivalent.

\section{Characterizing of minimal space-like surfaces in $\RR^4_1$ by means of $\Phi$ }

Let $\M$ be a surface in $\RR^4_1$, parameterized by isothermal
coordinates and let $\Phi$ be the function defined by
\eqref{Phi_def}.

Differentiating \eqref{Phi2}, we find:
\begin{equation}\label{Phi.dPhi_dbt}
\Phi\cdot\frac{\partial\Phi}{\partial\bar t}=0 \; .
\end{equation}
In view of \eqref{dPhi_dbt=dbPhi_dt} the function
$\ds\frac{\partial\Phi}{\partial\bar t}$ is real and
\eqref{Phi.dPhi_dbt} after a complex conjugation implies that:
\begin{equation}\label{bPhi.dPhi_dbt}
\bar\Phi\cdot\frac{\partial\Phi}{\partial\bar t}=0 \; .
\end{equation}
Since $\Phi$ and $\bar\Phi$ form a basis for $T_{C}(M)$, then
equalities \eqref{Phi.dPhi_dbt} and \eqref{bPhi.dPhi_dbt} mean
that $\ds\frac{\partial\Phi}{\partial\bar t}$ is orthogonal to
$T(\M)$ and therefore
\begin{equation}\label{dPhi_dbt_inN}
\frac{\partial\Phi}{\partial\bar t} \in N(\M) \; .
\end{equation}
In view of \eqref{dPhi_dbt} we find successively:
\begin{equation*}
\begin{array}{rl}\ds
\frac{\partial\Phi}{\partial\bar t}\!\! &=
\ds\left(\frac{\partial\Phi}{\partial\bar t}\right)^\bot=
\frac{1}{2}(\Delta \x )^\bot= \frac{1}{2}(\x_{uu}+\x_{vv} )^\bot=
\frac{1}{2}(\nabla_{\x_u} \x_u + \nabla_{\x_v} \x_v )^\bot\\[2.5ex]
&=\ds\frac{1}{2}(\sigma(\x_u,\x_u) + \sigma(\x_v,\x_v) )= E\;
\frac{1}{2}(\sigma(\vX_1,\vX_1) + \sigma(\vX_2,\vX_2) ) = E\vH\; .
\end{array}
\end{equation*}
Finally we have:
\begin{equation}\label{dPhi_dbt-Delta_x-EH}
\frac{\partial\Phi}{\partial\bar t}=\frac{1}{2}\Delta \x = E\vH\;
.
\end{equation}

Equality \eqref {dPhi_dbt-Delta_x-EH} implies immediately the
following statement.
\begin{thm}
Let $(\M,\x)$ be a space-like surface in $\RR^4_1$ parameterized by
isothermal coordinates $(u,v)\in \D$ and $\Phi(t)$ be the
complex-valued vector function in $\D$ defined by:
\[\Phi(t)=2\frac{\partial\x}{\partial t}=\x_u-\ii\x_v, \quad t=u+\ii . v,\]
Then the following conditions are equivalent:
\begin{enumerate}
    \item The function $\Phi(t)$ is holomorphic $\left( \ds\frac{\partial\Phi}{\partial\bar t}= 0 \right)$;%\\%[1ex]
    \item The function $\x (u,v)$ is harmonic $(\Delta \x = 0)$;
    \item $(\M,\x)$ is minimal space-like surface in $\RR^4_1$ $(\vH=0)$.
\end{enumerate}
\end{thm}

Let $(\M,\x)$ be a minimal space-like surface. Then the harmonic
conjugate function $\y$ to the function $\x$ is determined by the
Cauchy-Riemann equations
\[ \y_u=-\x_v; \quad \y_v=\x_u \ .\]
Let us introduce the function $\Psi$ by the equality
\[\Psi=\x+\ii\y .\]
The function $\Psi$ is holomorphic and $\x$, $\Phi$ are expressed
by $\Psi$ in the following way:
\[\x=\Re\Psi ; \qquad \Phi=\x_u-\ii\x_v=\x_u+\ii\y_u=\frac{\partial\Psi}{\partial u}=\Psi' .\]

Since $\ds\frac{\partial\Phi}{\partial\bar t}=0$, then
$\ds\frac{\partial\Phi}{\partial t}=\Phi'$ and
\begin{equation}\label{sigma_22}
\sigma(\vX_2,\vX_2)=-\sigma(\vX_1,\vX_1).
\end{equation}
Therefore
$$
\sigma(\x_v,\x_v)=E\sigma(\vX_2,\vX_2)=-E\sigma(\vX_1,\vX_1)=-\sigma(\x_u,\x_u).
$$

Then formulas \eqref{dPhi_dt} and \eqref{dPhi_dt_bot} get the
following form:
\begin{equation}\label{PhiPr}
\Phi^\prime=\frac{\partial\Phi}{\partial
u}=\x_{uu}-\ii\x_{uv};\quad \Phi^{\prime \bot}=\x_{uu}^\bot -\ii
x_{uv}^\bot =\sigma (\x_u,\x_u)-\ii\sigma (\x_u,\x_v).
\end{equation}

Formulas \eqref{sigma_uu_vv_uv} become correspondingly
\begin{equation}\label{sigma_uu_uv}
\begin{array}{l}
\sigma(\x_u,\x_u)=\ \ \,\Re (\Phi^{\prime \bot})=\ \
\,\ds\frac{1}{2}(\Phi^{\prime \bot}+\overline{\Phi^{\prime
\bot}})=
\ \ \,\ds\frac{1}{2}(\Phi^{\prime \bot}+{\overline{\Phi^\prime}}^\bot)\\[4mm]

\sigma(\x_v,\x_v)=     -\Re (\Phi^{\prime
\bot})=-\ds\frac{1}{2}(\Phi^{\prime \bot}+\overline{\Phi^{\prime
\bot}})=
-\ds\frac{1}{2}(\Phi^{\prime \bot}+{\overline{\Phi^\prime}}^\bot)\\[4mm]

\sigma(\x_u,\x_v)=-\Im (\Phi^{\prime \bot})=
\ds\frac{-1}{2\ii}(\Phi^{\prime \bot}-\overline{\Phi^{\prime
\bot}})= \ \ \!\,\ds\frac{\ii}{2}(\Phi^{\prime
\bot}-{\overline{\Phi^\prime}}^\bot).
\end{array}
\end{equation}

\section{Expressions for $K$ and $\varkappa$ of a minimal space-like surface by means of $\Phi$}

Let $(\M,\x)$ be a minimal space-like surface in $\RR^4_1$
parameterized by isothermal parameters. Choose a pair $\n_1$ and
$\n_2$ of orthonormal vector functions in $N(\M)$ of $\M$, such
that $\n^2_1=1$, $\n^2_2=-1$ and the quadruple
$(\vX_1,\vX_2,\n_1,\n_2)$ is right oriented in $\RR^4_1$.

For a given normal vector $\n$ we denote by $A_{\n}$ the
Weingarten operator in $T(\M)$. This operator and $\sigma$ are
related by the equality $A_{\n}\vX\cdot \vY=\sigma (\vX,\vY)\cdot
\n $. The condition $\vH=0$ implies that for any $\n$ $\trace
A_{\n}=0 $. Then the operators $A_{\n_1}$ and $A_{\n_2}$ have the
following representation
\begin{equation}\label{A1A2_nu_lambda_rho_mu}
A_{\n_1}= \left(
\begin{array}{rr}
\nu      &  \lambda\\
\lambda  & -\nu
\end{array}
\right); \qquad A_{\n_2}= \left(
\begin{array}{rr}
\rho &  \mu\\
\mu  & -\rho
\end{array}
\right)
\end{equation}
and the components of $\sigma$ are as follows:
\begin{equation}\label{sigma_nu_lambda_rho_mu}
\begin{array}{l}
\sigma (\vX_1,\vX_1)=(\sigma (\vX_1,\vX_1)\cdot \n_1)\n_1-(\sigma
(\vX_1,\vX_1)\cdot \n_2)\n_2=
\nu \n_1 - \rho \n_2\\
\sigma (\vX_1,\vX_2)=(\sigma (\vX_1,\vX_2)\cdot \n_1)\n_1-(\sigma
(\vX_1,\vX_2)\cdot \n_2)\n_2=
\lambda \n_1 - \mu \n_2\\
\sigma (\vX_2,\vX_2)=-\sigma (\vX_1,\vX_1)=-\nu \n_1 + \rho \n_2.
\end{array}
\end{equation}

Denote by $R$ the curvature tensor of ${\M}$. Then the Gauss
equation and \eqref{sigma_22} give:
\begin{equation}\label{K}
\begin{array}{rl}
K &=R(\vX_1,\vX_2)\vX_2\cdot \vX_1\\
  &=\sigma(\vX_1,\vX_1)\sigma(\vX_2,\vX_2)-\sigma^2(\vX_1,\vX_2)\\
  &=-\sigma^2(\vX_1,\vX_1)-\sigma^2 (\vX_1,\vX_2)\;.
\end{array}
\end{equation}
Now \eqref{sigma_nu_lambda_rho_mu} and \eqref{K} imply

\begin{equation}\label{K_nu_lambda_rho_mu}
K=-(\nu^2-\rho^2)-(\lambda^2-\mu^2)=\det(A_{\n_1})-\det(A_{\n_2})\;
.
\end{equation}
\medskip
On the other hand we get from \eqref{PhiPr}:
\begin{equation}\label{PhiPr_X}
\Phi^{\prime \bot}=\sigma ({\sqrt E}{\vX_1},{\sqrt
E}{\vX_1})-\ii\sigma ({\sqrt E}{\vX_1},{\sqrt E}{\vX_2})= E(\sigma
(\vX_1,\vX_1)-\ii\sigma (\vX_1,\vX_2)) \; .
\end{equation}
Calculating the norm $\|\Phi^{\prime \bot}\|$, we find
$$
\begin{array}{rl}{\|\Phi^{\prime \bot}\|}^2=\Phi^{\prime \bot}\cdot\overline{\Phi^{\prime \bot}}
&=E(\sigma (\vX_1,\vX_1)-\ii\sigma (\vX_1,\vX_2))E(\sigma (\vX_1,\vX_1)+\ii\sigma (\vX_1,\vX_2))\\[2mm]
&=E^2(\sigma^2(\vX_1,\vX_1)+\sigma^2 (\vX_1,\vX_2)).\end{array}
$$
Taking into account the last equality and \eqref{EG} we have:
\begin{equation}\label{s2+s2}
\sigma^2(\vX_1,\vX_1)+\sigma^2 (\vX_1,\vX_2)=\frac{{\|\Phi^{\prime
\bot}\|}^2}{E^2}=\frac{4{\|\Phi^{\prime \bot}\|}^2}{\|\Phi\|^4}\;.
\end{equation}
Now \eqref{s2+s2} and  \eqref{K} imply that
\begin{equation}\label{K_Phi}
K= \ds\frac{-4{\|\Phi^{\prime \bot}\|}^2}{\|\Phi\|^4}.
\end{equation}

 In the last formula we can represent $\|\Phi^{\prime \bot}\|^2$ in a different way. Note that $\Phi^2=0$
means that $\Phi$ and $\bar\Phi$ are orthogonal with respect to the
Hermitian dot product in ${\CC}^4$. In view of \eqref{Phi_def} and
\eqref{xuxv} it follows that they form an orthogonal basis of the
complexified tangential plane of ${\M}$. Therefore the tangential
projection of $\Phi^\prime$ is as follows:
$$\Phi^{\prime\top}
=\ds\frac{\Phi^{\prime\top}\cdot\bar
\Phi}{\|\Phi\|^2}\Phi+\ds\frac{\Phi^{\prime\top}\cdot \Phi}{\|\bar
\Phi\|^2}\bar \Phi =\ds\frac{\Phi' \cdot \bar
\Phi}{\|\Phi\|^2}\Phi + \ds\frac{\Phi' \cdot \Phi}{\|\bar
\Phi\|^2}\bar \Phi. $$ Differentiating $\Phi^2=0$, we find
$\Phi\cdot\Phi^\prime=0$. Thus the projection of $\Phi'$ has the
form:
\begin{equation}\label{Phipn}
\Phi^{\prime\top}= \ds\frac{\Phi' \cdot \bar
\Phi}{\|\Phi\|^2}\Phi; \quad\quad
\Phi^{\prime\bot}=\Phi'-\Phi^{\prime\top}= \Phi'-\ds\frac{\Phi'
\cdot \bar \Phi}{\|\Phi\|^2}\Phi.
\end{equation}
Using the second equality of \eqref{Phipn} by means of complex
conjugation we get:
\begin{equation*}%\label{mPhipn2}
\begin{array}{rl}
{\|\Phi^{\prime\bot}\|}^2&=\Phi^{\prime\bot}\cdot\overline{\Phi^{\prime\bot}}=
\left(\Phi'-\ds\frac{\Phi' \cdot \bar \Phi}{\|\Phi\|^2}\Phi\right)\left(\bar{\Phi'}-\ds\frac{\bar{\Phi'} \cdot \Phi}{\|\Phi\|^2}\bar \Phi \right)\\[4mm]
                      &=\Phi'\cdot\bar{\Phi'}-\ds\frac{\bar{\Phi'} \cdot \Phi}{\|\Phi\|^2}\Phi'\cdot\bar \Phi -
                                              \ds\frac{\Phi' \cdot \bar \Phi}{\|\Phi\|^2}\Phi\cdot\bar{\Phi'} +
                                                \ds\frac{(\Phi' \cdot \bar \Phi)(\bar{\Phi'} \cdot \Phi)}{\|\Phi\|^4}\Phi\cdot\bar \Phi\\[4mm]
                                            &=\|\Phi'\|^2 - \ds\frac{|\bar{\Phi'} \cdot \Phi|^2}{\|\Phi\|^2} - \ds\frac{|\Phi' \cdot \bar \Phi|^2}{\|\Phi\|^2}+
                                              \ds\frac{|\Phi' \cdot \bar \Phi|^2}{\|\Phi\|^4}\|\Phi\|^2 = \|\Phi'\|^2 - \ds\frac{|\bar{\Phi'} \cdot \Phi|^2}{\|\Phi\|^2}\\[4mm]
                                            &=\ds\frac{\|\Phi\|^2\|\Phi'\|^2-|\bar \Phi \cdot \Phi'|^2}{\|\Phi\|^2}\ .
\end{array}
\end{equation*}
Since the norm of the bi-vector $\Phi\wedge\Phi'$ is given by:
\[\|\Phi\wedge\Phi'\|^2=\|\Phi\|^2\|\Phi'\|^2-|\bar \Phi \cdot \Phi'|^2,\]
then we have
\begin{equation*}
\|\Phi^{\prime\bot}\|^2=\ds\frac{\|\Phi\|^2\|\Phi'\|^2-|\bar \Phi
\cdot \Phi'|^2}{\|\Phi\|^2}=
\ds\frac{\|\Phi\wedge\Phi'\|^2}{\|\Phi\|^2}\
\end{equation*}
and
\begin{equation}\label{K_Phi_bv}
K= \ds\frac{-4{\|\Phi^{\prime \bot}\|}^2}{\|\Phi\|^4}=
   \ds\frac{-4\|\Phi\wedge\Phi'\|^2}{\|\Phi\|^6}\ .
\end{equation}
\medskip
In order to obtain formulas for the normal curvature $\varkappa$,
let us denote by $R^N$ the curvature tensor of the normal
connection of ${\M}$. The Ricci equation and
\eqref{A1A2_nu_lambda_rho_mu} imply
\begin{equation}\label{kappa}
\begin{array}{rl}
\varkappa &= R^N(\vX_1,\vX_2,\n_1,\n_2)=R^N(\vX_1,\vX_2)\n_2\cdot
\n_1 \\ \;
          &=[A_{\n_2},A_{\n_1}] \vX_1\cdot \vX_2=A_{\n_1}\vX_1\cdot A_{\n_2}\vX_2-A_{\n_2}\vX_1\cdot A_{\n_1}\vX_2\\
                    &=(\nu\vX_1 + \lambda\vX_2)\cdot (\mu\vX_1 - \rho\vX_2)-(\rho\vX_1 + \mu\vX_2)\cdot (\lambda\vX_1 - \nu\vX_2)\\
                    &=\nu\mu - \nu\rho + \lambda\mu - \lambda\rho  - (\rho\lambda - \rho\nu + \mu\lambda - \mu\nu) \\
                    &=2\nu\mu - 2\rho\lambda \; .
\end{array}
\end{equation}
We denote by $\det (\va,\vb,\vc,\vd)$ the determinant of the
vectors $\va$, $\vb$, $\vc$ and $\vd$, with respect to the
standard basis in ${\CC}^4$. Taking into account
\eqref{sigma_nu_lambda_rho_mu}, we have
$$
\begin{array}{rl}
\det (\x_u,\x_v,\sigma(\x_u,\x_u),\sigma(\x_u,\x_v))&=
\det (\sqrt E \vX_1,\sqrt E \vX_2,\sigma(\sqrt E \vX_1,\sqrt E \vX_1),\sigma(\sqrt E \vX_1,\sqrt E \vX_2))\\
                                              &=E^3 \det (\vX_1,\vX_2,\sigma(\vX_1,\vX_1),\sigma(\vX_1,\vX_2))\\
                                              &=E^3 \det (\vX_1,\vX_2,\nu\n_1 - \rho\n_2,\lambda\n_1 - \mu\n_2)\\
                                                                                            &=-E^3 \det (\vX_1,\vX_2,\nu\n_1,\mu\n_2)-E^3 \det (\vX_1,\vX_2,\rho\n_2,\lambda\n_1)\\
                                              &=E^3(-\nu\mu + \rho\lambda)\det (\vX_1,\vX_2,\n_1,\n_2)=E^3(-\nu\mu + \rho\lambda).$$
\end{array}
$$
From the last equation it follows that
\begin{equation}\label{numu2}
-\nu\mu + \rho\lambda = \ds \frac{1}{E^3}\det
(\x_u,\x_v,\sigma(\x_u,\x_u),\sigma(\x_u,\x_v)).
\end{equation}
Replacing $\x_u$ and $\x_v$ by \eqref{xuxv} we find
\begin{equation}\label{det1}
\begin{array}{l}
\ds \det (\x_u,\x_v,\sigma(\x_u,\x_u),\sigma(\x_u,\x_v))=\frac{\ii}{4}\det (\Phi+\bar\Phi,\Phi-\bar\Phi,\sigma(\x_u,\x_u),\sigma(\x_u,\x_v))\\
\ds =\frac{\ii}{4}\det (\Phi,-\bar\Phi,\sigma(\x_u,\x_u),\sigma(\x_u,\x_v))+ \frac{\ii}{4}\det (\bar\Phi,\Phi,\sigma(\x_u,\x_u),\sigma(\x_u,\x_v))\\
\ds =-\frac{\ii}{2}\det
(\Phi,\bar\Phi,\sigma(\x_u,\x_u),\sigma(\x_u,\x_v)).
\end{array}
\end{equation}
Similarly, using \eqref{sigma_uu_uv}, we get:
\begin{equation}\label{det2}
\begin{array}{l}
\ds \det
(\Phi,\bar\Phi,\sigma(\x_u,\x_u),\sigma(\x_u,\x_v))=-\frac{\ii}{2}\det
(\Phi,\bar\Phi,\Phi^{\prime \bot},{\overline{\Phi^\prime}}^\bot).
\end{array}
\end{equation}
Now \eqref{det2} and \eqref{det1} give
\begin{equation*}
\begin{array}{l}
\ds \det (\x_u,\x_v,\sigma(\x_u,\x_u),\sigma(\x_u,\x_v))=-\frac{1}{4}\det (\Phi,\bar\Phi,\Phi^{\prime \bot},{\overline{\Phi^\prime}}^\bot)\\
\ds =-\frac{1}{4}\det (\Phi,\bar\Phi,\Phi^\prime -\Phi^{\prime
\top},\overline{\Phi^\prime}-{\overline{\Phi^\prime}}^\top).
\end{array}
\end{equation*}
Hence
\begin{equation}\label{det3}
\begin{array}{l}
\ds \det (\x_u,\x_v,\sigma(\x_u,\x_u),\sigma(\x_u,\x_v))=
-\frac{1}{4}\det
(\Phi,\bar\Phi,\Phi^\prime,\overline{\Phi^\prime}).
\end{array}
\end{equation}

In view of \eqref{kappa}, \eqref{numu2} and \eqref{det3} we have

\begin{equation*}
\begin{array}{l}
\ds \varkappa = 2\nu\mu - 2\rho\lambda = \frac{-2}{E^3}\det
(\x_u,\x_v,\sigma(\x_u,\x_u),\sigma(\x_u,\x_v))=
\frac{1}{2E^3}\det
(\Phi,\bar\Phi,\Phi^\prime,\overline{\Phi^\prime}).
\end{array}
\end{equation*}
Using \eqref{EG} we find:
\begin{equation}\label{kappa2}
\begin{array}{l}
\ds \varkappa = \frac{4}{\|\Phi\|^6}\det
(\Phi,\bar\Phi,\Phi^\prime,\overline{\Phi^\prime}).
\end{array}
\end{equation}

For any minimal space-like surface $(\M,\x)$ in $\RR^4_1$,
parameterized by isothermal coordinates the Gauss curvature $K$
and the normal curvature $\varkappa$ are given by the formulas:

\begin{equation}\label{K_kappa_nu_lambda_rho_mu}
K= -\nu^2-\lambda^2+\rho^2+\mu^2, \quad\quad \varkappa = 2\nu\mu -
2\rho\lambda \; ;
\end{equation}
\begin{equation}\label{K_kappa_Phi}
K= \ds\frac{-4{\|\Phi^{\prime \bot}\|}^2}{\|\Phi\|^4}
=\ds\frac{-4\|\Phi\wedge\Phi'\|^2}{\|\Phi\|^6}, \quad\quad
\varkappa = \ds\frac{4}{\|\Phi\|^6}\det
(\Phi,\bar\Phi,\Phi^\prime,\overline{\Phi^\prime}).
\end{equation}

\section{Existence of canonical coordinates on a minimal space-like surface}

Let $\M$ be a minimal space-like surface in $\RR^4_1$.
\begin{defn}\label{def_izr}
A point $p\in\M$ is said to be degenerate if the set \\
$\{\sigma(\vX,\vY);\ \vX \in T_p(\M),\vY \in T_p(\M) \}$, is
contained into one of the light-like one-dimensional subspaces of
$N_p(\M)$.
\end{defn}
Let $(\M ,\ \x=\Re\Psi)$ be a minimal space-like surface in
$\RR^4_1$ parameterized by isothermal coordinates $(u,v)$.

\begin{thm}\label{izr_Phi}
A point $p\in \M$ is degenerate if and only if ${\Phi^{\prime
\bot}}^2=0$.
\end{thm}
Proof. Let us consider again equality \eqref{PhiPr_X}. Squaring
both sides of the equality, we find
\begin{equation}\label{PhiPr_bot2}
\begin{array}{rl}
{\Phi^{\prime \bot}}^2 \!\!\! &=E^2(\sigma^2(\vX_1,\vX_1)-\ii\, 2\sigma (\vX_1,\vX_1)\sigma (\vX_1,\vX_2)-\sigma^2 (\vX_1,\vX_2))\\
                       &=E^2(\sigma^2(\vX_1,\vX_1)-\sigma^2 (\vX_1,\vX_2))-\ii\, 2E^2 \sigma (\vX_1,\vX_1)\sigma (\vX_1,\vX_2).
\end{array}
\end{equation}

The last equality implies the following equivalence:
\begin{equation}\label{izr}
\begin{array}{l}
\sigma (\vX_1,\vX_1)\bot \: \sigma (\vX_1,\vX_2)\\[2mm]
\sigma^2(\vX_1,\vX_1)=\sigma^2 (\vX_1,\vX_2)
\end{array}
\quad \Leftrightarrow \quad {\Phi^{\prime \bot}}^2=0.
\end{equation}

Assuming that the point into consideration is degenerate, then it
follows that $\sigma^2(\vX_1,\vX_1)=0$, $\sigma^2(\vX_1,\vX_2)=0$
and $\sigma (\vX_1,\vX_1)\sigma (\vX_1,\vX_2)=0$. Now \eqref{izr}
implies that ${\Phi^{\prime \bot}}^2=0$.

Let ${\Phi^{\prime \bot}}^2=0$. We have to prove that the vectors
$\sigma(\vX_1,\vX_1)$ and $\sigma(\vX_1,\vX_2)$ lie in one and the
same light-like one-dimensional subspace of $N_p(\M)$. If we
assume that $\sigma^2(\vX_1,\vX_1)>0$, then it follows that
$\sigma^2(\vX_1,\vX_2)>0$ and $\sigma (\vX_1,\vX_1)\bot \: \sigma
(\vX_1,\vX_2)$, which is a contradiction.

Similarly, assuming that $\sigma^2(\vX_1,\vX_1)<0$, we obtain the
metric on $N_p(\M)$ is negative definite, which is a
contradiction.

Thus $\sigma^2(\vX_1,\vX_1)=0$, $\sigma^2(\vX_1,\vX_2)=0$ and
$\sigma (\vX_1,\vX_1)\bot \: \sigma (\vX_1,\vX_2)$. Hence the
vectors $\sigma(\vX_1,\vX_1)$ and $\sigma(\vX_1,\vX_2)$ are
light-like and lie in one and the same one-dimensional subspace of
$N_p(\M)$. \prfend

Next we prove that ${\Phi^{\prime \bot}}^2$ is a holomorphic
function of $t$. In general it doesn't follow that the projection
$\Phi^{\prime \bot}$ is a holomorphic function, but we shall prove
that ${\Phi^{\prime \bot}}^2={\Phi^\prime}^2$.

In order to prove the last equality, we square the second equality
in \eqref{Phipn} and get:
$$
{\Phi^{\prime\bot}}^2 = {\Phi'}^2-2\Phi'\ds\frac{\Phi' \cdot \bar
\Phi}{\|\Phi\|^2}\Phi+\left(\ds\frac{\Phi' \cdot \bar
\Phi}{\|\Phi\|^2}\right)^2 \Phi^2.
$$
Taking into account equalities $\Phi^2=0$ and
$\Phi\cdot\Phi^\prime=0$, we find:
\begin{equation}\label{PhiPr_bot2=PhiPr2}
{\Phi^{\prime \bot}}^2={\Phi^\prime}^2 \; .
\end{equation}

Thus we obtained that any degenerate point of $\M$ is a zero of
the holomorphic function ${\Phi^\prime}^2$. This implies
immediately the following characterization of the set of
degenerate points of a minimal space-like surface:
\begin{thm}
If $\M$ is a connected minimal space-like surface in $\RR^4_1$,
then: either it consists of degenerate points or the set of the
degenerate points is countable without any limit points.
\end{thm}

Further in this section we consider minimal space-like surfaces in
$\RR^4_1$ without degenerate points.

We give the following definitions:
\begin{defn}\label{canonical1}
The isothermal coordinates $(u, v)$ on a minimal space-like
surface are said to be \emph{canonical of the first type} if
\begin{equation}\label{can1}
\begin{array}{l}
\sigma (\vX_1,\vX_1)\bot \: \sigma (\vX_1,\vX_2),\\[2mm]
E^2(\sigma^2(\vX_1,\vX_1)-\sigma^2 (\vX_1,\vX_2))=1.
\end{array}
\end{equation}
\end{defn}

Because of \eqref{PhiPr_bot2} the isothermal parameters $(u,v)$
are canonical of the first type if and only if
\begin{equation}\label{can1phi}
{\Phi^\prime}^2={\Phi^{\prime \bot}}^2=1.
\end{equation}

\begin{defn}\label{canonical2}
The isothermal coordinates $(u, v)$ on a minimal space-like
surface are said to be \emph{canonical of the second type} if
\begin{equation}\label{can2}
\begin{array}{l}
\sigma (\vX_1,\vX_1)\bot \: \sigma (\vX_1,\vX_2),\\[2mm]
E^2(\sigma^2(\vX_1,\vX_1)-\sigma^2 (\vX_1,\vX_2))=-1.
\end{array}
\end{equation}
\end{defn}

The isothermal coordinates $(u,v)$ are canonical of the second
type if and only if:
\begin{equation}\label{can2phi}
{\Phi^\prime}^2={\Phi^{\prime \bot}}^2=-1.
\end{equation}

\begin{thm}\label{can_coord_exist}
Any minimal space-like surface in $\RR^4_1$, free of degenerate
points, admits locally canonical coordinates of both types.
\end{thm}

Proof. For arbitrary isothermal coordinates $(u,v)$ on the
surface, denote $t=u+v\ii$. Let us consider the change $t=t(\tilde
t\:)$ of the complex variable $t$ by the new complex variable
$\tilde t\:$. We shall find the conditions under which the new
variable determines canonical coordinates. First, the new
coordinates have to be isothermal, i.e. $t=t(\tilde t\,)$ is a
conformal map in ${\CC}$. Therefore the function $t(\tilde t\,)$
is either holomorphic or antiholomorphic. The case of an
antiholomorphic function is reduced to the case of a holomorphic
function by means of the additional change ${\tilde t}=\bar s$.

It is enough to consider only the case of a holomorphic function
$t(\tilde t\,)$. Let $\tilde\Psi$ be the holomorphic function
representing ${\M}$ with respect to the new coordinates and
$\tilde\Phi$ be its derivative. Then we have
\begin{equation}\label{tildPhi}
\tilde\Phi = \tilde\Psi'_{\tilde t} = \Psi'_t t' = \Phi t'.
\end{equation}
The derivative of $\tilde\Phi$ with respect to $\tilde t$ is given
by $\tilde\Phi'_{\tilde t}=\Phi'_t t'^2+\Phi t''$. Since $\Phi$ is
tangent to the surface ${\M}$, then $\Phi^\bot=0$ and consequently
\begin{equation}\label{tildPhiPr2}
\begin{array}{lll}
\tilde\Phi_{\tilde t}'^\bot &=& (\Phi'_t t^{\prime \, 2}+\Phi t'')^\bot = \Phi_t'^\bot t^{\prime \, 2},\\
\left.\tilde\Phi_{\tilde t}^{\prime\bot}\right.^2 &=&
{\Phi_t^{\prime \bot}}^2 t'^4.
\end{array}
\end{equation}
According to \eqref{can1} and \eqref{can2} the new complex
variable $\tilde t$ determines canonical coordinates if
$\left.\tilde\Phi_{\tilde t}^{\prime\bot}\right.^2=\pm 1$. If
${\Phi_t^{\prime \bot}}^2=0$, then by virtue of \eqref{tildPhiPr2}
it follows that $\left.\tilde\Phi_{\tilde
t}^{\prime\bot}\right.^2=0$. The last condition means that the
point is degenerate, which is impossible. Hence ${\Phi^{\prime
\bot}}^2 \neq 0$. Then the function $\tilde t$ determines
canonical coordinates if and only if ${\Phi_t^{\prime \bot}}^2
t'^4=\pm 1$, i.e. $t(\tilde t\,)$ satisfies the following ordinary
complex first order differential equation:
\begin{equation}\label{eqcan}
\sqrt[4]{\pm{\Phi_t^{\prime \bot}}^2}\:dt = d{\tilde t}.
\end{equation}
Integrating \eqref{eqcan} and taking into account that the left
side of the equality is holomorphic, we obtain $\tilde t$ as a
holomorphic function of $t$. The condition ${\Phi_t^{\prime
\bot}}^2 \neq 0$ means that ${\tilde t}'\neq 0$ and the
correspondence between ${\tilde t}$ and $t$ is one-to-one.
Consequently ${\tilde t}$ determines isothermal coordinates
satisfying the condition $\left.\tilde\Phi_{\tilde
t}^{\prime\bot}\right.^2=\pm 1$, which implies that they are
canonical. \prfend

Next we consider the question of uniqueness of canonical
coordinates. Suppose that $t$ and $\tilde t$ are canonical of one
and the same type. Then $t = t(\tilde t\,)$ is either holomorphic
or antiholomorphic. According to \eqref{can1} and \eqref{can2}
equality \eqref{tildPhiPr2} implies that
$$
\pm 1 = \left.\tilde\Phi_{\tilde
t}^{\prime\bot}\right.^2={\Phi_t^{\prime \bot}}^2 t'^4 = \pm 1
t'^4 = \pm t'^4.
$$
Therefore $t'^4 = 1$ and $t' = \pm 1;\ \pm \ii$. We get from here
that $t$ and $\tilde t$ are related by one of the following
equalities: $t=\pm\tilde t+c;\ \pm \ii\tilde t+c$, where
$c=\text{const}$.

The anti-holomorphic case is reduced to the holomorphic one by the
change ${\tilde t}=\bar s$ and we get: $t=\pm\bar{\tilde t}+c;\ \pm \ii\bar{\tilde t}+c$.

Thus we obtain eight possible relations between $t$ and $\tilde t$. Under the natural initial
condition $c=0$, these relations mean that:

\emph{The canonical coordinates of one and the same type are
unique up to a direction and numbering of the coordinate lines.}

Finally, we consider the relations between canonical coordinates
of different type. Let $t=u+v\ii$ be canonical coordinates of the
first type and introduce new coordinates by means of the formula
$t = e^{\frac{\pi \ii}{4}}\tilde t$. Then $t'^4=-1$,
$\left.\tilde\Phi_{\tilde t}^{\prime\bot}\right.^2=-1$, and
consequently $\tilde t$ determines canonical coordinates of the
second type. This construction shows that the canonical
coordinates of both types are obtained from each other by a
rotation of the angle $\frac{\pi}{4}$ in the coordinate plane
$(u,v)$.

Let $(\M,\x)$ be a minimal space-like surface in $\RR^4_1$ free of
degenerate points, parameterized by canonical coordinates of the
first type. We can precise the choice of the orthonormal pair
{$\n_1, \n_2$ in $N(\M)$. Since $\sigma (\vX_1,\vX_1)\bot \:
\sigma (\vX_1,\vX_2)$, then we can choose $\n_1$ and $\n_2$ to be
collinear with $\sigma (\vX_1,\vX_1)$ and $\sigma (\vX_1,\vX_2)$.
More precisely, if at a point we have $\sigma (\vX_1,\vX_1)\neq
0$, then we choose $\n_1$ with the same direction as $\sigma
(\vX_1,\vX_1)$, and $\n_2$ so that the quadruple
$(\vX_1,\vX_2,\n_1,\n_2)$ is a positive oriented basis in
$\RR^4_1$. Then $\n_2$ is collinear with $\sigma (\vX_1,\vX_2)$.
Under these conditions formulas \eqref{sigma_nu_lambda_rho_mu} get
the form:
\begin{equation}\label{sigma_nu_mu}
\begin{array}{l}
\sigma (\vX_1,\vX_1)=\phantom{-} \nu \, \n_1, \\
\sigma (\vX_1,\vX_2)=         -  \mu \, \n_2, \\
\sigma (\vX_2,\vX_2)=         -  \nu \, \n_1;
\end{array}
\qquad \nu>0 \; .
\end{equation}
Therefore we have $\lambda=0$ and $\rho=0$ and formulas
\eqref{A1A2_nu_lambda_rho_mu} become as follows:
\begin{equation}\label{A1A2_nu_mu}
A_{\n_1}= \left(
\begin{array}{rr}
\nu  &  0\\
0    & -\nu
\end{array}
\right); \qquad A_{\n_2}= \left(
\begin{array}{rr}
0    & \mu\\
\mu  & 0
\end{array}
\right).
\end{equation}

If at a given non-degenerate point $\sigma (\vX_1,\vX_1)=0$, then
$\sigma (\vX_1,\vX_2)\neq 0$. In this case we can choose first
$\n_2$ collinear with the same direction with $-\sigma
(\vX_1,\vX_2)$, and then $\n_1$ so that the quadruple
$(\vX_1,\vX_2,\n_1,\n_2)$ forms a positive oriented basis in
$\RR^4_1$.

The functions $\nu$ and $\mu$ also satisfy the following
relations:
\begin{equation}\label{numu1}
\begin{array}{lr}
\nu^2 \!\!  &= \phantom{-} \sigma^2 (\vX_1,\vX_1),\\
\mu^2 \!\!  &=            -\sigma^2 (\vX_1,\vX_2).
\end{array}
\end{equation}

The functions $\nu$ and $\mu$ are a pair of scalar invariants of a
minimal space-like surface, free of degenerate points. These
invariants completely determine the second fundamental form via
\eqref{sigma_nu_mu}. The second condition in \eqref{can1} implies
that the first fundamental form is completely determined by the
formula:
\begin{equation}\label{E_nu_mu}
E=\frac{1}{\sqrt{\nu^2+\mu^2}} \; .
\end{equation}

The relations between the pairs $(\nu ,\mu)$ and $(K ,\varkappa)$
are as follows:
\begin{equation}\label{K_kappa_nu_mu}
K= -\nu^2+\mu^2\; , \quad\quad \varkappa = 2\nu\mu \; .
\end{equation}

\begin{equation}\label{nu_mu_K_kappa}
\mu^2=\frac{\sqrt{K^2+\varkappa^2}+K}{2}\,, \qquad
\nu^2=\frac{\sqrt{K^2+\varkappa^2}-K}{2}\,.
\end{equation}
Using the above formulas we can characterize the degenerate points
of ${\M}$ in terms of $K$ and $\varkappa$.

\begin{thm}
Let $\M$ be a minimal space-like surface with Gaussian curvature
$K$ and normal curvature $\varkappa$. A point $p\in \M$ is
degenerate if and only if $K=0$ and $\varkappa =0$.
\end{thm}
Proof. If $p$ is not a degenerate point in $\M$, then we can
introduce canonical coordinates of the first type in a
neighborhood of $p$. Formulas \eqref{numu1} imply that at least
one of $(\nu ,\mu )$ is different from $0$. Applying
\eqref{K_kappa_nu_mu} we obtain that at least one of $(K
,\varkappa )$ is also different from 0.

If $p$ is a degenerate point, then $\sigma (\vX_1,\vX_1)$ and
$\sigma (\vX_1,\vX_2)$ are lightlike. Then \eqref{K} implies that
$K=0$. Further it follows that $\sigma (\vX_1,\vX_1)$ and $\sigma
(\vX_1,\vX_2)$ are collinear. Therefore the determinant of the
four vectors $\x_u$,$\x_v$,$\sigma(\x_u,\x_u)$ and
$\sigma(\x_u,\x_v)$ is zero. Hence, in view of \eqref{det3} and
\eqref{kappa2} it follows that $\varkappa =0$. \prfend

Finally we add some formulas for $\nu$, $\mu$ and $\varkappa$ in
canonical coordinates of the first type. Equalities \eqref{can1}
and \eqref{EG} imply that
\begin{equation}\label{s2-s2}
\sigma^2(\vX_1,\vX_1)-\sigma^2
(\vX_1,\vX_2)=\frac{1}{E^2}=\frac{4}{\|\Phi\|^4}\;.
\end{equation}

By virtue of \eqref{s2-s2}, \eqref{s2+s2} and \eqref{numu1} we
find
\begin{equation}\label{nu2+-mu2}
\nu^2 \!\!  = \ds \frac{2(1+{\|\Phi^{\prime
\bot}\|}^2)}{\|\Phi\|^4}\,, \qquad \mu^2 \!\!  = \ds
\frac{2(1-{\|\Phi^{\prime \bot}\|}^2)}{\|\Phi\|^4}\,.
\end{equation}
Hence
\begin{equation}\label{modkappa}
|\varkappa \:|=|2\nu \mu |= \ds \frac{4\sqrt{{1-\|\Phi^{\prime
\bot}\|}^4}}{\|\Phi\|^4}.
\end{equation}

\section{General Weierstrass representations for minimal space-like surfaces.}

In this section we give several types of general Weierstrass representations for
minimal space-like surfaces in $\RR^4_1$. In $\RR^4$ such formulas were considered in
\cite{HO}, \cite{G-K-2}. In $\RR^4_1$ general Weierstrass representations were used
in \cite{AA-VJ_DegGM}, \cite{AA-VJ_PrescGM}.

Let $(\M,\x)$:\; $\x=\Re\Psi$ be a minimal space-like surface in
$\RR^4_1$, parameterized by isothermal coordinates, and let
$\Phi=\Psi^\prime$. If $\Phi=(\phi_1,\phi_2,\phi_3,\phi_4)$, then
the condition $\Phi^2=0$ is equivalent to that the coordinates are
\begin{equation}\label{phi2coord}
\phi_1^2+\phi_2^2+\phi_3^2-\phi_4^2=0
\end{equation}

The relation \eqref{phi2coord} ca be 'parameterized' in different
ways by means of three holomorphic functions.

First we shall find a representation of $\Phi$ by means of
trigonometric functions.

We write \eqref{phi2coord} in the following forms:
$$\phi_1^2+\phi_2^2=-\phi_3^2+\phi_4^2, \quad
  \phi_1^2+\phi_3^2=-\phi_2^2+\phi_4^2, \quad
  \phi_2^2+\phi_3^2=-\phi_1^2+\phi_4^2.$$
At least one of these three quantities $\phi_1^2+\phi_2^2$,
$\phi_1^2+\phi_3^2$ and $\phi_2^2+\phi_3^2$ is different from
zero. (The opposite leads to $\phi_1^2+\phi_2^2+\phi_3^2=0$ and
$\phi_1^2=\phi_2^2=\phi_3^2=\phi_4^2=0$, which contradicts to the
condition ${\M}$ is regular.) Without loss of generality we can
assume that $\phi_1^2+\phi_2^2\neq 0$, which means that there
exists a holomorphic function $f \neq 0$, such that:
\begin{equation}\label{f}
f^2=\phi_1^2+\phi_2^2=-\phi_3^2+\phi_4^2.
\end{equation}
The last equality is equivalent to the next one:
\begin{equation}\label{f2}
\left(\frac{\phi_1}{f}\right)^2+\left(\frac{\phi_2}{f}\right)^2=
\left(\frac{\phi_3}{\ii
f}\right)^2+\left(\frac{\phi_4}{f}\right)^2=1.
\end{equation}
Hence, there exist holomorphic functions $h_1$ and $h_2$, such
that:
$$
\frac{\phi_1}{f}=\cos h_1, \quad  \frac{\phi_2}{f}=\sin h_1, \quad
\frac{\phi_3}{\ii f}=\cos h_2, \quad \frac{\phi_4}{f}=\sin h_2.
$$
Thus we found the following representation of the function $\Phi$:
\begin{equation}\label{W1}
\Phi: \quad
\begin{array}{rlr}
\phi_1 &=& f\cos h_1,\\
\phi_2 &=& f\,\sin h_1,\\
\phi_3 &=& \ii f\cos h_2,\\
\phi_4 &=& f\,\sin h_2;\\
\end{array}
\qquad f \neq 0.
\end{equation}
Next we have to express the condition $\|\Phi\|^2>0$ in terms of
the triple $(f,h_1,h_2)$. Equality \eqref{W1} implies that
\begin{equation}\label{mPhi2_fh1h2}
\begin{array}{rl}
\|\Phi\|^2=\Phi\bar\Phi &= |f|^2(\cos h_1 \cos \bar h_1 + \sin h_1 \sin \bar h_1\\
                               &+\cos h_2 \cos \bar h_2 - \sin h_2 \sin \bar h_2)\\
                &= |f|^2(\cos  (h_1 - \bar h_1) + \cos  (h_2 + \bar h_2))\\
                            &= |f|^2(\cos  (2\ii \Im h_1) + \cos  (2\Re h_2))\\
                                &= |f|^2(\cosh (2\Im h_1 )    + \cos (2\Re h_2 )).
\end{array}
\end{equation}
Since $\cosh (2\Im h_1 )\ge 1 \ge |\cos (2\Re h_2 )|$, then it
follows from \eqref{mPhi2_fh1h2} that $\|\Phi\|^2\ge 0$. The
equality is equivalent to $\cosh (2\Im h_1 )=1$ and $\cos (2\Re
h_2 )=-1$, i.e. $\Im h_1=0$ and $\Re h_2=\frac{\pi}{2}+k\pi;\ k\in
{\mathbb Z}$. Thus we obtained that the triple $(f,h_1,h_2)$ in
the representation \eqref{W1} satisfies the conditions
\begin{equation}\label{W1_cond}
f\neq 0; \quad \Im h_1\neq 0\ \text{or}\ \Re h_2\neq
\frac{\pi}{2}+k\pi;\ k\in {\mathbb Z}.
\end{equation}

Hence, any minimal space-like surface ${\M}$ in $\RR^4_1$,
parameterized by isothermal coordinates, admits Weierstrass
representation of the type \eqref{W1}, where the triple $(f,h_1,h_2)$ satisfies
the conditions \eqref{W1_cond}.

Conversely, any triple $(f,h_1,h_2)$ of holomorphic functions,
defined in a domain in ${\CC}$ and satisfying the conditions
\eqref{W1_cond}, determines by \eqref{W1} a holomorphic
${\CC}^4$-valued function $\Phi$. It follows from \eqref{W1_cond}
that $\|\Phi\|^2>0$. By direct computations we get $\Phi^2=0$.
Then the surface ${\M}:\; \x=\Re (\Psi)$, where $\Psi$ is
determined by the equality $\Psi^\prime = \Phi$, is a minimal
space-like surface in $\RR^4_1$, parameterized by isothermal coordinates.

So, we proved the following statement. \vskip 1mm

\emph{Any triple of holomorphic functions $(f,h_1,h_2)$ satisfying
\eqref{W1_cond}, generates by means of formulas \eqref{W1} a
minimal space-like surface in $\RR^4_1$.}

\vskip 1mm Finally let us establish to what extent the function
$\Phi$ determines the functions $(f,h_1,h_2)$.

Suppose that one and the same function $\Phi$ is represented by
two different triples $(f,h_1,h_2)$ and $(\hat f,\hat h_1,\hat
h_2)$. Then \eqref{f} and \eqref{W1} imply the following relations
between both triples:
\begin{equation*}
\begin{array}{ll}
\hat f   \!\! &=\: f\\
\hat h_1 \!\! &=\: h_1 + 2k_1\pi\\
\hat h_2 \!\! &=\: h_2 + 2k_2\pi
\end{array}\quad
\text{or}\quad
\begin{array}{ll}
\hat f   \!\! &=\: -f\\
\hat h_1 \!\! &=\: h_1 + (2k_1+1)\pi\\
\hat h_2 \!\! &=\: h_2 + (2k_2+1)\pi
\end{array};\quad \ \
\begin{array}{l}
k_1\in {\mathbb Z}\\
k_2\in {\mathbb Z}
\end{array}
\end{equation*}
\vspace {2mm}

With the aid of different substitutions in \eqref{W1} we can
obtain other Weierstrass representations for minimal space-like
surfaces in $\RR^4_1$.

In order to obtain Weierstrass representation by means of
hyperbolic functions, we make the following substitution in
\eqref{W1}:
$$
f\rightarrow \ii f; \quad h_1 \rightarrow -\ii h_1; \quad h_2
\rightarrow \pi +\ii h_2.
$$
Thus we obtain the following Weierstrass representation:
\begin{equation}\label{W2}
\Phi: \quad
\begin{array}{rlr}
\phi_1 &=& \ii f\cosh h_1,\\
\phi_2 &=&  f\,\sinh h_1,\\
\phi_3 &=&  f\cosh h_2,\\
\phi_4 &=&  f\,\sinh h_2.\\
\end{array}
\end{equation}

Taking into account \eqref{W1_cond}, it follows that the functions
$(f,h_1,h_2)$ satisfy the conditions:
\begin{equation}\label{W2_cond}
f\neq 0; \quad \Re h_1\neq 0\ \text{or}\ \Im h_2\neq
\frac{\pi}{2}+k\pi;\ k\in {\mathbb Z}.
\end{equation}

Further, let us change the functions $h_1$ and $h_2$ in
$\eqref{W2}$ by $w_1$ and $w_2$ in the following way:
\begin{equation}\label{w1w2}
\begin{array}{l}
w_1=h_1+h_2\\
w_2=h_1-h_2
\end{array}
\end{equation}
Then we obtain the following representation of the surface:
\begin{equation}\label{W5}
\Phi: \quad
\begin{array}{rlr}
\phi_1 &=& \ii f \cosh \ds\frac{w_1+w_2}{2}\\[4mm]
\phi_2 &=&  f \sinh \ds\frac{w_1+w_2}{2}\\[4mm]
\phi_3 &=&  f \cosh \ds\frac{w_1-w_2}{2}\\[4mm]
\phi_4 &=&  f \sinh \ds\frac{w_1-w_2}{2}\\
\end{array}
\end{equation}

It follows from \eqref{W2_cond} that $(f,w_1,w_2)$ satisfy the conditions:
\begin{equation}\label{W5_cond}
f\neq 0; \quad \Re (w_1 + w_2) \neq 0 \ \text{or}\ \Im (w_1 - w_2) \neq (2k+1)\pi;\ k\in {\mathbb Z}.
\end{equation}
The last conditions can be written in the form:
\begin{equation*}
f\neq 0; \quad \Re (w_1 + \bar w_2) \neq 0 \ \text{or}\ \Im (w_1 + \bar w_2) \neq (2k+1)\pi;\ k\in {\mathbb Z}.
\end{equation*}
Thus we obtained the following more simple form for the conditions \eqref{W5_cond}:
\begin{equation}\label{W5_cond-2}
f\neq 0; \quad w_1 + \bar w_2 \neq (2k+1)\pi\ii;\ k\in {\mathbb Z}.
\end{equation}

Next we introduce the functions $g_1$ and $g_2$ by the equalities:
\begin{equation}\label{g}
g_1=e^{w_1}; \quad g_2=e^{w_2}.
\end{equation}
Using these functions, we obtain from \eqref{W5} the Weierstrass
representation, which is the analogue of the classical Weierstrass
representation for minimal surfaces in $\RR^3$. Consequently we
calculate the coordinate functions:
$$\phi_1 =\frac{\ii f}{2} (e^{\frac{w_1+w_2}{2}}+e^{-\frac{w_1+w_2}{2}})=
          \frac{\ii f}{2\sqrt{g_1 g_2}}(g_1 g_2+1),$$
$$\phi_2 = \frac{f}{2} (e^{\frac{w_1+w_2}{2}}-e^{-\frac{w_1+w_2}{2}})=
         \frac{f}{2\sqrt{g_1 g_2}}(g_1 g_2-1).$$
$$\phi_3 = \frac{f}{2} (e^{\frac{w_1-w_2}{2}}+e^{-\frac{w_1-w_2}{2}})=
           \frac{f}{2\sqrt{g_1 g_2}}(g_1+g_2).$$
$$\phi_4 = \ds\frac{f}{2} (e^{\frac{w_1-w_2}{2}}-e^{-\frac{w_1-w_2}{2}})=
         \ds\frac{f}{2\sqrt{g_1 g_2}}(g_1 - g_2).$$
\noindent
In the last equalities we make the substitution
\begin{equation}\label{ff}
f \rightarrow f2\sqrt{g_1 g_2}
\end{equation}
and obtain the following 'polynomial' Weierstrass representation:
\begin{equation}\label{W6}
\Phi: \quad
\begin{array}{rll}
\phi_1 &=& \ii f(g_1 g_2+1),\\
\phi_2 &=& \ f(g_1 g_2-1),\\
\phi_3 &=& \ f(g_1+g_2),\\
\phi_4 &=& \ f(g_1-g_2).\\
\end{array}
\end{equation}
Now we shall determine the conditions which satisfy the functions
$(f, g_1, g_2)$. It follows from \eqref{g} that
the condition $w_1 + \bar w_2 \neq (2k+1)\pi\ii;\ k\in {\mathbb Z}$ is
equivalent to the condition
$e^{w_1 + \bar w_2} \neq e^{(2k+1)\pi\ii}=-1;\ k\in {\mathbb Z}$,
which gives
$g_1 \bar g_2 \neq -1$.
Therefore we obtained from \eqref{W5_cond-2} the following conditions:
\begin{equation}\label{W6_cond-2}
f\neq 0; \quad g_1 \bar g_2 \neq -1.
\end{equation}

Conversely, if $(f,g_1,g_2)$ are three holomorphic functions defined in a domain in ${\CC}$
and satisfying \eqref{W6_cond-2}, then formulas \eqref{W6} determine a holomorphic function
$\Phi$ with values in ${\CC}^4$. Equalities \eqref{W6_cond-2} imply that $\|\Phi\|^2>0$. By
direct computations we get from \eqref{W6} equality \eqref{phi2coord}, which is $\Phi^2=0$.
If we determine the function $\Psi$ by the equality $\Psi^\prime = \Phi$ and define
${\M}:\; \x=\Re (\Psi)$, then ${\M}$ is a minimal space-like surface in $\RR^4_1$,
parameterized by isothermal coordinates.

Thus we obtained:

{\it Any three holomorphic functions $(f,g_1,g_2)$ satisfying \eqref{W6_cond-2}, generates via \eqref{W6}
a minimal space-like surface in $\RR^4_1$.}
\begin{rem}
We obtained the representation \eqref{W6} using \eqref{g}, which implies that the functions
$g_1$ and $g_2$ are different from zero at any point. This follows from the fact that we chose
$\phi_1^2+\phi_2^2\neq 0$. If any of the functions $(g_1, g_2)$ is zero at a fixed point,
then it follows directly from \eqref{W6} $\Phi^2=0$ and $\|\Phi\|^2=2|f|$. This means that
\eqref{W6} again determine a minimal space-like surface in $\RR^4_1$. Therefore there is
no need to add new conditions for $g_1$ and $g_2$ other than these from \eqref{W6_cond-2}.
\end{rem}

In the end we show that the functions $(f,g_1,g_2)$ can be
expressed by the components of the vector function $\Phi$.
Directly from \eqref{W6} we get:
$$
\begin{array}{l}
\ii\phi_1+\phi_2=-f(g_1 g_2+1)+f(g_1 g_2-1)=-2f,\\
\phi_3+\phi_4=f(g_1+g_2)+f(g_1-g_2)=2fg_1,\\
\phi_3-\phi_4=f(g_1+g_2)-f(g_1-g_2)=2fg_2.
\end{array}
$$
Hence, the functions $f$, $g_1$ and $g_2$ are expressed as
follows:
\begin{equation}\label{fg1g2}
f=-\ds\frac{1}{2}(\ii\phi_1+\phi_2), \quad
g_1=-\ds\frac{\phi_3+\phi_4}{\ii\phi_1+\phi_2}\,, \quad
g_2=-\ds\frac{\phi_3-\phi_4}{\ii\phi_1+\phi_2}\,.
\end{equation}
\section{Some formulas, related to Weierstrass representations}
In this section we use the Weierstrass representation \eqref{W2}
for minimal space-like surfaces by means of hyperbolic functions.
Using the functions $f, h_1, h_2$, respectively $f, w_1, w_2$, we
obtain some formulas, which we use further.

First we introduce some subsidiary functions and denotations.

The holomorphic vector function $\va$ is defined by the equality:
\begin{equation}\label{adef}
\va=\ds\frac{\Phi}{f}.
\end{equation}

Next we introduce the following denotations:
\begin{equation}\label{alpha_befta}
\alpha=\Re (h_1), \quad \beta=\Im (h_2).
\end{equation}
These functions determine the function $\theta$, given by:
\begin{equation}\label{theta_def}
\theta = \Re h_1 + \ii\Im h_2 = \alpha + \ii\beta .
\end{equation}
The function $\theta$ is a complex harmonic function, which in
general is not holomorphic.

Under these denotations applying the Cauchy-Riemann equations, we
have:
\begin{equation}\label{hpr_albet}
\begin{array}{l}
h'_1=\Re(h_1)'_u+\ii\Im(h_1)'_u=\Re(h_1)'_u-\ii\Re(h_1)'_v=\alpha'_u-\ii\alpha'_v\;,\\
h'_2=\Re(h_2)'_u+\ii\Im(h_2)'_u=\Im(h_2)'_v+\ii\Im(h_2)'_u=\beta'_v+\ii\beta'_u\;.
\end{array}
\end{equation}
For $w'_1$ and $w'_2$ we find, respectively:
\begin{equation}\label{wpr_albet}
\begin{array}{l}
w'_1=(\alpha'_u+\beta'_v)-\ii (\alpha'_v-\beta'_u)\;,\\
w'_2=(\alpha'_u-\beta'_v)-\ii (\alpha'_v+\beta'_u)\;.
\end{array}
\end{equation}
Using \eqref{W2} and \eqref{adef}, we get the following formulas
for $\va$, $\bar \va$, $\va'$ and $\bar {\va'}$:
\begin{equation}\label{a}
\begin{array}{l}
\va=(\ \ \, \ii\cosh h_1,\sinh h_1,\cosh h_2,\sinh h_2)\\
\bar \va=(-\ii\cosh \bar h_1,\sinh \bar h_1,\cosh \bar h_2,\sinh \bar h_2)\\
\va'\! =(\ \ \, \ii h'_1 \sinh h_1,h'_1 \cosh h_1,h'_2 \sinh h_2,h'_2 \cosh h_2)\\
\bar{\va'}\! =(-\ii\bar{h'_1} \sinh \bar h_1,\bar{h'_1} \cosh \bar
h_1,\bar{h'_2} \sinh \bar h_2,\bar{h'_2} \cosh \bar h_2)
\end{array}
\end{equation}

Further we find the scalar products between the functions $\va$,
$\bar \va$, $\va'$ and $\bar {\va'}$. Differentiating the equality
$\va^2=0$, we have:
\begin{equation}\label{a1}
\va^2=\va \va'=\bar \va^2=\bar \va \bar {\va'}=0
\end{equation}
Taking scalar multiplications in \eqref{a}, we also obtain:
\begin{equation}\label{a2}
\begin{array}{rl}
\|\va\|^2=\va\bar \va &= \cosh h_1 \cosh \bar h_1 + \sinh h_1 \sinh \bar h_1 + \cosh h_2 \cosh \bar h_2 - \sinh h_2 \sinh \bar h_2\\
                &= \cosh (h_1 + \bar h_1) + \cosh (h_2 - \bar h_2)\\
                                &= \cosh (2\Re h_1 ) + \cosh (2\ii\Im h_2 )\\
                                &= 2\cosh(\Re h_1+ \ii\Im h_2)\cosh( \Re h_1- \ii\Im h_2)\\
                                &= 2\cosh(\theta)\cosh(\bar\theta)=2|\cosh(\theta)|^2;
\end{array}
\end{equation}
\begin{equation}\label{a3}
\begin{array}{rl}
\va\bar {\va'} &= \bar{h'_1}\cosh h_1 \sinh \bar h_1 +
\bar{h'_1}\sinh h_1 \cosh \bar h_1 +
              \bar{h'_2}\cosh h_2 \sinh \bar h_2 - \bar{h'_2}\sinh h_2 \cosh \bar h_2\\
           &= \bar{h'_1}\sinh (h_1 + \bar h_1) - \bar{h'_2}\sinh (h_2 - \bar h_2)\\
                     &= \bar{h'_1}\sinh (2\Re h_1 ) - \bar{h'_2}\sinh (2\ii\Im h_2 );
\end{array}
\end{equation}
\begin{equation}\label{a4}
 \bar \va \va' = \overline{\va\bar {\va'}} = h'_1 \sinh (2\Re h_1 ) + h'_2 \sinh (2\ii\Im h_2 );
\end{equation}
\begin{equation}\label{a5}
\begin{array}{rl}
\va'^2 &= -h'^2_1\sinh^2 h_1 + h'^2_1\cosh^2 h_1 + h'^2_2\sinh^2 h_2 - h'^2_2\cosh^2 h_2\\
     &= h'^2_1-h'^2_2 = w'_1 w'_2;
\end{array}
\end{equation}
\begin{equation}\label{a6}
\begin{array}{rl}
\|\va'\|^2 = \va'\bar {\va'} &= |h'_1|^2 \sinh h_1 \sinh \bar h_1 + |h'_1|^2 \cosh h_1 \cosh \bar h_1\\
                       &+\  |h'_2|^2 \sinh h_2 \sinh \bar h_2 - |h'_2|^2 \cosh h_2 \cosh \bar h_2\\
                       &= |h'_1|^2\cosh (h_1 + \bar h_1) - |h'_2|^2\cosh (h_2 - \bar h_2)\\
                                 &= |h'_1|^2\cosh (2\Re h_1 ) - |h'_2|^2\cosh (2\ii\Im h_2 ).
\end{array}
\end{equation}
Further we obtain formulas for $\va^{\prime\bot}$,
${\va^{\prime\bot}}^2$ and ${\|\va^{\prime\bot}\|}^2$ expressed by
means of $h_1, h_2$ and $w_1, w_2$, respectively. For
$\va^{\prime\bot}$ we have
$\va^{\prime\bot}=\va'-\va^{\prime\top}$. The equality $\va^2=0$
means that the vectors $\va$ and $\bar \va$ are orthogonal with
respect to the Hermitian dot product in ${\CC}^4$. Therefore the
tangential vector $\va^{\prime\top}$ is decomposed as follows:
$$\va^{\prime\top}=\ds\frac{\va^{\prime\top}\cdot\bar \va}{\|\va\|^2}\va+\ds\frac{\va^{\prime\top}\cdot \va}{\|\bar \va\|^2}\bar \va
                =\ds\frac{\va' \cdot \bar \va}{\|\va\|^2}\va + \ds\frac{\va' \cdot \va}{\|\bar \va\|^2}\bar \va. $$
Equality \eqref{a1} implies that $\va' \cdot \va = 0$. Thus we
obtained:
\begin{equation}\label{apn}
\va^{\prime\top}= \ds\frac{\va' \cdot \bar \va}{\|\va\|^2}\va;
\quad\quad \va^{\prime\bot}=\va'-\va^{\prime\top}=
\va'-\ds\frac{\va' \cdot \bar \va}{\|\va\|^2}\va.
\end{equation}
Taking square in both sides of \eqref{apn}, we get:
$$
{\va^{\prime\bot}}^2 = {\va'}^2-2\va'\ds\frac{\va' \cdot \bar
\va}{\|\va\|^2}\va+\left(\ds\frac{\va' \cdot \bar
\va}{\|\va\|^2}\right)^2 \va^2.
$$
Taking again into account \eqref{a1}, we have $\va' \cdot \va = 0$
and $\va^2=0$. Consequently ${\va^{\prime\bot}}^2 = {\va'}^2$. Now
by virtue of \eqref{a5} we find
\begin{equation}\label{apn2}
{\va^{\prime\bot}}^2 = {\va'}^2 = {h'_1}^2-{h'_2}^2 = w'_1 w'_2
\end{equation}

Using \eqref{apn} and applying complex conjugation, we calculate
${\|\va^{\prime\bot}\|}^2$:
\begin{equation}\label{mapn2}
\begin{array}{rl}
{\|\va^{\prime\bot}\|}^2&=\va^{\prime\bot}\cdot\overline{\va^{\prime\bot}}=
\left(\va'-\ds\frac{\va' \cdot \bar \va}{\|\va\|^2}\va\right)\left(\bar{\va'}-\ds\frac{\bar{\va'} \cdot \va}{\|\va\|^2}\bar \va \right)\\[4mm]
                      &=\va'\cdot\bar{\va'}-\ds\frac{\bar{\va'} \cdot \va}{\|\va\|^2}\va'\cdot\bar \va -
                                              \ds\frac{\va' \cdot \bar \va}{\|\va\|^2}\va\cdot\bar{\va'} +
                                                \ds\frac{(\va' \cdot \bar \va)(\bar{\va'} \cdot \va)}{\|\va\|^4}\va\cdot\bar \va\\[4mm]
                                            &=\|\va'\|^2 - \ds\frac{|\bar{\va'} \cdot \va|^2}{\|\va\|^2} - \ds\frac{|\va' \cdot \bar \va|^2}{\|\va\|^2}+
                                              \ds\frac{|\va' \cdot \bar \va|^2}{\|\va\|^4}\|\va\|^2 = \|\va'\|^2 - \ds\frac{|\bar{\va'} \cdot \va|^2}{\|\va\|^2}\\[4mm]
                                            &=\ds\frac{\|\va\|^2\|\va'\|^2-|\bar \va \cdot \va'|^2}{\|\va\|^2}
\end{array}
\end{equation}

Let us denote the numerator in \eqref{mapn2} by $k_1$. Applying
equalities \eqref{a2}, \eqref{a4} and \eqref{a6} we find:
\begin{equation}\label{k1}
\begin{array}{rl}
k_1 &= \|\va\|^2\|\va'\|^2-|\bar \va \cdot \va'|^2\\
    &= (|h'_1|^2-|h'_2|^2)(1+\cosh (2\Re h_1 )\cos (2\Im h_2 ))\\
        &+ 2\Im(\bar h'_1 h'_2)\sinh (2\Re h_1 )\sin (2\Im h_2 )\\
        &= ({\alpha'_u}^2+{\alpha'_v}^2-{\beta'_u}^2-{\beta'_v}^2)(1+\cosh (2\alpha )\cos (2\beta ))\\
        &+ 2(\alpha'_u\beta'_u+\alpha'_v\beta'_v)\sinh (2\alpha )\sin (2\beta )\;.
\end{array}
\end{equation}
Denote the determinant of the vectors $\va$, $\bar \va$, $\va'$
and $\bar {\va'}$ by $k_2$. Applying formulas \eqref{a}, we find:
\begin{equation}\label{k2}
\begin{array}{rl}
k_2 &=  \det(\va,\bar \va , \va' , \bar {\va'})\\
    &= -2\Im(\bar h'_1 h'_2)(1+\cosh (2\Re h_1 )\cos (2\Im h_2 ))\\
        &+ (|h'_1|^2-|h'_2|^2)\sinh (2\Re h_1 )\sin (2\Im h_2 )\\
        &= -2(\alpha'_u\beta'_u+\alpha'_v\beta'_v)(1+\cosh (2\alpha )\cos (2\beta ))\\
        &+ ({\alpha'_u}^2+{\alpha'_v}^2-{\beta'_u}^2-{\beta'_v}^2)\sinh (2\alpha )\sin (2\beta )\;.
\end{array}
\end{equation}

Next we simplify the expressions for $k_1$ and $k_2$ calculating
the complex quantity $-k_1+\ii k_2$:
\begin{equation}\label{-k1+ik2}
\begin{array}{rl}
-k_1+\ii k_2 &=
-\,({\alpha'_u}^2+{\alpha'_v}^2-{\beta'_u}^2-{\beta'_v}^2+2\ii
(\alpha'_u\beta'_u+\alpha'_v\beta'_v))
                 (1+\cosh (2\alpha )\cos (2\beta ))\\
        &\phantom{=}+(2\ii (\alpha'_u\beta'_u+\alpha'_v\beta'_v)+{\alpha'_u}^2+{\alpha'_v}^2-{\beta'_u}^2-{\beta'_v}^2)
                             \sinh (2\alpha )\sinh (2\ii\beta )\\
                       &= -\,((\alpha'_u+\ii\beta'_u)^2+(\alpha'_v+\ii\beta'_v)^2)\\
                &\phantom{=} \phantom{-}\ \, (1+\cosh (2\alpha )\cosh (2\ii\beta )-\sinh (2\alpha )\sinh (2\ii\beta ))\\
                         &= -2((\alpha'_u+\ii\beta'_u)^2+(\alpha'_v+\ii\beta'_v)^2)\cosh^2(\alpha-\ii\beta)\;.
\end{array}
\end{equation}
Using the function $\theta$, defined by \eqref{theta_def}, we
obtain another form of $-k_1+\ii k_2$:
\begin{equation}\label{-k1+ik2_theta}
-k_1+\ii k_2 =
-2({\theta'_u}^2+{\theta'_v}^2)\cosh^2(\bar\theta)\;.
\end{equation}
Further we express $-k_1+\ii k_2$ in terms of $w_1$ and $w_2$. For
the first factor in \eqref{-k1+ik2_theta} we have:
\begin{equation*}
\begin{array}{rl}
{\theta'_u}^2+{\theta'_v}^2&=(\alpha'_u+\ii\beta'_u)^2+(\alpha'_v+\ii\beta'_v)^2\\
                    &=(\alpha'_u+\ii\beta'_u+\ii (\alpha'_v+\ii\beta'_v))(\alpha'_u+\ii\beta'_u-\ii (\alpha'_v+\ii\beta'_v))\\
                                        &=(\alpha'_u-\beta'_v+\ii (\alpha'_v+\beta'_u))(\alpha'_u+\beta'_v-\ii (\alpha'_v-\beta'_u))\;.
\end{array}
\end{equation*}
Comparing the last formula with \eqref{wpr_albet}, we get:
\begin{equation*}
{\theta'_u}^2+{\theta'_v}^2=(\Re w'_2 + \ii (-\Im w'_2))(\Re w'_1
+ \ii \Im w'_1)\;.
\end{equation*}
The above formulas imply that:
\begin{equation}\label{thetapr_wpr}
{\theta'_u}^2+{\theta'_v}^2=w'_1 \bar {w'_2}\;.
\end{equation}
In order to find the second factor in \eqref{-k1+ik2_theta}, first
we find $\theta$:
\begin{equation*}
\begin{array}{rl}
\theta &= \alpha+\ii\beta = \Re h_1 + \ii\Im h_2\\
       &= \frac{1}{2}(h_1 + \bar h_1) + \ii \frac{1}{2\ii}(h_2 - \bar h_2)\\[0.4ex]
             &= \frac{1}{2}(h_1 + h_2) + \frac{1}{2}(\bar h_1 - \bar h_2)\;.
\end{array}
\end{equation*}
Taking into account the above equality and \eqref{w1w2}, we find:
\begin{equation}\label{theta_w12}
\theta =\frac{w_1 + \bar w_2}{2}\;.
\end{equation}
Consequently
\begin{equation*}
\begin{array}{rl}
\cosh(\theta) &=\frac{1}{2}(e^{\frac{w_1 + \bar
w_2}{2}}+e^{-\frac{w_1 + \bar w_2}{2}})
                   =\frac{1}{2}e^{-\frac{w_1 + \bar w_2}{2}}(e^{w_1 + \bar w_2}+1)\\
                                    &=\frac{1}{2}e^{-\frac{w_1}{2}}e^{-\frac{\bar w_2}{2}}(1+e^{w_1}e^{\bar w_2})\;.

\end{array}
\end{equation*}
Finally we have
\begin{equation}\label{cosh^2theta_w12}
\cosh^2\theta =\frac{1}{4}e^{- w_1}e^{- \bar w_2}(1+e^{w_1}e^{\bar
w_2})^2\;.
\end{equation}

Now we replace \eqref{thetapr_wpr} and \eqref{cosh^2theta_w12}
into \eqref{-k1+ik2_theta} and obtain:
\begin{equation}\label{-k1+ik2_w12}
-k_1+\ii k_2 = -\frac{1}{2}w'_1 \bar {w'_2}e^{-\bar w_1}e^{-
w_2}(1+e^{\bar w_1}e^{w_2})^2\;.
\end{equation}

\section{Canonical Weierstrass representation for minimal space-like surfaces of general type}

In this section we introduce canonical Weierstrass representations for minimal
space-like surfaces of general type in $\RR^4_1$. Weierstrass representations with
respect to canonical coordinates were obtained in \cite{G-2} for $\RR^3_1$ and in
\cite{G-K-2} for $\RR^4$.

\begin{defn}\label{def_gentype}
A minimal space-like surface in $\RR^4_1$ is said to be of
\emph{general type} if it is free of degenerate points in the
sense  of \ref{def_izr}.
\end{defn}

Let the minimal space-like surface ${\M}$ of general type be
parameterized by canonical coordinates of the first type. Consider
the Weierstrass representation \eqref{W2} by means of hyperbolic
functions. The condition \eqref{can1} leads to a relation between
the three functions $f$, $h_1$ and $h_2$. In order to obtain this
relation, we express the condition ${\Phi^{\prime \bot}}^2=1$ via
$f$, $h_1$ and $h_2$. By virtue of \eqref{adef} we have
$\Phi=f\va$ and therefore $\Phi'=f'\va+f\va'$. Since the vector
$\va$ is tangential to ${\M}$, then we get:
\begin{equation}\label{Phipna}
\Phi^{\prime \bot}=(f'\va+f\va')^\bot = f\va^{\prime \bot}; \quad
\quad {\Phi^{\prime \bot}}^2=f^2 {\va^{\prime \bot}}^2.
\end{equation}
Because of \eqref{apn2} we have ${\va^{\prime\bot}}^2 =
{h'_1}^2-{h'_2}^2$ and consequently ${\Phi^{\prime \bot}}^2=f^2
({h'_1}^2-{h'_2}^2)$. Taking into account the last equality and
\eqref{PhiPr_bot2=PhiPr2}, we obtain that the minimal space-like
surface ${\M}$ given by \eqref{W2} is parameterized by canonical
coordinates of the first type if and only if:
\begin{equation}\label{can1h}
{\Phi'}^2=f^2 ({h'_1}^2-{h'_2}^2)=1\;.
\end{equation}
The last formula and \eqref{W2} imply the following statement.
\begin{thm}
Any minimal space-like surface ${\M}$ of general type,
parameterized by canonical coordinates of the first type, has the
following Weierstrass representation:
\begin{equation}\label{Wcanh}
\Phi: \quad
\begin{array}{rlr}
\phi_1 &=& \ii\ds\frac{\cosh h_1}{\sqrt{{h'_1}^2 - {h'_2}^2}}\,,\\[8mm]
\phi_2 &=&  \ds\frac{\sinh h_1}{\sqrt{{h'_1}^2 - {h'_2}^2}}\,,\\[8mm]
\phi_3 &=&  \ds\frac{\cosh h_2}{\sqrt{{h'_1}^2 - {h'_2}^2}}\,,\\[8mm]
\phi_4 &=&  \ds\frac{\sinh h_2}{\sqrt{{h'_1}^2 - {h'_2}^2}}\,,\\
\end{array}
\end{equation}
where $(h_1,h_2)$ are holomorphic functions satisfying the
conditions:
\begin{equation}\label{Wcanh_cond}
{h'_1}^2 \neq {h'_2}^2; \quad \Re h_1\neq 0\ \text{or}\ \Im
h_2\neq \frac{\pi}{2}+k\pi;\ k\in {\mathbb Z}.
\end{equation}

Conversely, if $(h_1,h_2)$ is a pair of holomorphic functions
satisfying the conditions \eqref{Wcanh_cond}, then formulas
\eqref{Wcanh} give a minimal space-like surface of general
type, parameterized by canonical coordinates of the first type.
\end{thm}
We call the representation of $\Phi$ in Theorem 9.2
\emph{canonical Weierstrass representation}.

Using the functions $w_1$ and $w_2$, given by \eqref{w1w2}, then
the condition \eqref{can1h} gets the form:
\begin{equation}\label{can1w}
{\Phi'}^2=f^2 w'_1 w'_2 = 1
\end{equation}
If we replace  $h_1$ and $h_2$ with $w_1$ and $w_2$ into
\eqref{Wcanh}, then we obtain the following canonical Weierstrass
representation for ${\M}$:
\begin{equation}\label{Wcanw}
\Phi: \quad
\begin{array}{rlr}
\phi_1 &=& \ds\frac{\ii}{\sqrt{w'_1 w'_2}} \cosh \ds\frac{w_1+w_2}{2}\,,\\[6mm]
\phi_2 &=& \ds\frac{1}{\sqrt{w'_1 w'_2}} \sinh \ds\frac{w_1+w_2}{2}\,,\\[6mm]
\phi_3 &=& \ds\frac{1}{\sqrt{w'_1 w'_2}} \cosh \ds\frac{w_1-w_2}{2}\,,\\[6mm]
\phi_4 &=& \ds\frac{1}{\sqrt{w'_1 w'_2}} \sinh \ds\frac{w_1-w_2}{2}\,.\\
\end{array}
\end{equation}
According to \eqref{W5_cond-2}, the functions $(w_1,w_2)$ satisfy
the conditions:
\begin{equation}\label{Wcanw_cond}
w'_1 w'_2 \neq 0; \quad w_1 + \bar w_2 \neq (2k+1)\pi\ii;\ k\in {\mathbb Z}.
\end{equation}

Conversely, if $(w_1,w_2)$ is a pair of holomorphic functions,
satisfying the conditions \eqref{Wcanw_cond}, then the formulas
\eqref{Wcanw} generate a minimal space-like surface of general
type, parameterized by canonical coordinates of the first type.

Finally, using the functions $g_1$ and $g_2$, given by \eqref{g},
we obtain a canonical Weierstrass representation of the type
\eqref{W6}. Differentiating \eqref{g}, we get:
\begin{equation}\label{gp}
g'_1=e^{w_1}w'_1=g_1 w'_1; \quad g'_2=e^{w_2}w'_2=g_2 w'_2.
\end{equation}
From here we have:
\begin{equation}\label{wp}
w'_1=\frac{g'_1}{g_1}\,, \quad w'_2=\frac{g'_2}{g_2}.
\end{equation}
Applying \eqref{ff} and \eqref{wp} to the condition \eqref{can1w},
we get $(f2\sqrt{g_1 g_2})^2 \ds\frac{g'_1}{g_1}
\ds\frac{g'_2}{g_2}=1$.

Consequently the isothermal coordinates are canonical of the first
type if and only if
\begin{equation}\label{can1g}
{\Phi'}^2=4f^2 g'_1 g'_2 = 1.
\end{equation}
Next we express $f$ from the last equality of \eqref{can1g} and
replace it into \eqref{W6}. Thus we obtain the following statement.

\begin{thm}\label{thm_Wcang}
Any minimal space-like surface ${\M}$ of general type,
parameterized by canonical coordinates of the first type, has the
following Weierstrass representation:
\begin{equation}\label{Wcang}
\Phi: \quad
\begin{array}{rll}
\phi_1 &=& \ds\frac{\ii}{2}\; \ds\frac {g_1 g_2+1}{\sqrt{g'_1 g'_2}}\,,\\[6mm]
\phi_2 &=& \ds\frac{1}{2}\; \ds\frac {g_1 g_2-1}{\sqrt{g'_1 g'_2}}\,,\\[6mm]
\phi_3 &=& \ds\frac{1}{2}\; \ds\frac {g_1 + g_2}{\sqrt{g'_1 g'_2}}\,,\\[6mm]
\phi_4 &=& \ds\frac{1}{2}\; \ds\frac {g_1 - g_2}{\sqrt{g'_1 g'_2}}\,.\\
\end{array}
\end{equation}
According to \eqref{W6_cond-2} the functions $(g_1,g_2)$ in this representation satisfy
the conditions:
\begin{equation}\label{Wcang_cond}
g'_1 g'_2 \neq 0; \quad g_1 \bar g_2 \neq -1.
\end{equation}

Conversely, if $(g_1,g_2)$ is a pair of holomorphic functions
satisfying the conditions \eqref{Wcang_cond}, then formulas
\eqref{Wcang} generate a minimal space-like surface of general
type, parameterized by canonical coordinates of the first type.
\end{thm}
The above canonical Weierstrass representation seems to be the
most useful and applicable representation.

\section{The first fundamental form and the curvatures $K, \varkappa$ in a general Weierstrass representation}

Let ${\M}$ be a minimal space-like surface in $\RR^4_1$,
parameterized by isothermal coordinates. First we consider the
case, when ${\M}$ is given by \eqref{W2}. In order to obtain a
formula for $E$, we use equalities \eqref{EG}, \eqref{adef} and
\eqref{a2}. Thus we get:
\begin{equation}\label{E_ftheta}
E
=\frac{1}{2}\|\Phi\|^2=\frac{1}{2}\|f\va\|^2=|f|^2|\cosh\theta|^2\;,
\end{equation}
where $\theta$ is the function \eqref{theta_def}.

Taking into account \eqref{theta_w12}, we find the following
formula for $E$ with respect to the representation \eqref{W5}:
\begin{equation}\label{E_fw1w2}
E=|f|^2\left|\cosh\frac{w_1 + \bar w_2}{2}\right|^2.
\end{equation}
In order to obtain a formula with respect to the representation
\eqref{W6}, we express $\cosh(\theta)$ by means of $g_j$,
$(j=1;2)$, given by \eqref{g}. Then \eqref{cosh^2theta_w12} gives
that
\begin{equation}\label{cosh^2theta_g12}
\cosh^2(\theta) =\frac{(1+g_1\bar g_2)^2}{4g_1\bar g_2}\;.
\end{equation}
Passing from the representation \eqref{W5} by means of $w_j$ to
the representation \eqref{W6} by means of $g_j$, $(j=1;2)$, as a
consequence of \eqref{ff} $|f|^2$ has to be replaced by:
\begin{equation}\label{mf2_mf2}
|f|^2 \rightarrow 4|f|^2|g_1 g_2|\;.
\end{equation}
Now applying \eqref{cosh^2theta_g12} and \eqref{mf2_mf2} to
\eqref{E_ftheta}, we find the following formula for the
coefficient $E$ in the representation \eqref{W6}:
\begin{equation}\label{E_fg1g2}
E=|f|^2|1+g_1\bar g_2|^2\;.
\end{equation}
%Thus we obtained that the condition \eqref{W6_cond} in the
%representation \eqref{W6} is equivalent to
%\begin{equation}\label{W6_cond_2}
%f\neq 0; \quad g_1\bar g_2 \neq -1\;.
%\end{equation}
%The condition \eqref{Wcang_cond} in the canonical representation
%\eqref{Wcang} is equivalent respectively to
%\begin{equation}\label{Wcang_cond_2}
%g'_1 g'_2 \neq 0; \quad g_1\bar g_2 \neq -1\;.
%\end{equation}

Further we find the corresponding formulas for $K$ and $\varkappa$.
In the formula \eqref{K_kappa_Phi} we replace $\Phi^{\prime \bot}$
by \eqref{Phipna} and get:
\[
K = \ds\frac{-4{\|\Phi^{\prime \bot}\|}^2}{\|\Phi\|^4} =
\ds\frac{-4{\|f\va^{\prime \bot}\|}^2}{\|f\va\|^4} =
\ds\frac{-4{|f|^2\|\va^{\prime \bot}\|}^2}{|f|^4\|\va\|^4} =
\ds\frac{-4{\|\va^{\prime \bot}\|}^2}{|f|^2\|\va\|^4}.
\]
Using \eqref{mapn2} and taking into account \eqref{k1}, we find:
\[
K = \ds\frac{-4(\|\va\|^2\|\va'\|^2-|\bar \va \cdot
\va'|^2)}{|f|^2\|\va\|^6} = \ds\frac{-4k_1}{|f|^2\|\va\|^6}.
\]

A similar formula for $\varkappa$ can be derived using the second
equality in \eqref{K_kappa_Phi}. We find consecutively:
\[
\ds \varkappa =  \ds\frac{4}{\|f\va\|^6}\det (f\va,\bar f \bar
\va,f'\va+f\va',\bar f' \bar \va + \bar f \bar {\va'})
              =  \ds\frac{4|f|^4}{|f|^6\|\va\|^6}\det (\va,\bar \va,\va',\bar {\va'}) = \ds\frac{4k_2}{|f|^2\|\va\|^6}.
\]
Thus we have:
\begin{equation}\label{Kkappa1}
K = \ds\frac{-4k_1}{|f|^2\|\va\|^6}; \quad \varkappa =
\ds\frac{4k_2}{|f|^2\|\va\|^6}\;.
\end{equation}

It is useful to unite $K$ and $\varkappa$ in one formula by the
complex quantity $K+\ii\varkappa$. It follows from \eqref{Kkappa1}
that:
\begin{equation}\label{K+ikappa1}
K+\ii\varkappa = \ds\frac{4(-k_1+\ii k_2)}{|f|^2\|\va\|^6}\;.
\end{equation}

Replacing $\|\va\|^2$ and $-k_1+\ii k_2$ respectively by
\eqref{a2} and \eqref{-k1+ik2_theta} we get:
\[
K+\ii\varkappa  =
\ds\frac{4(-2({\theta'_u}^2+{\theta'_v}^2)\cosh^2(\bar\theta))}{|f|^2
8|\cosh(\theta)|^6}\;.
\]
Finally we obtained:
\begin{equation}\label{K+ikappa_ftheta}
K+\ii\varkappa  = \ds\frac{-({\theta'_u}^2+{\theta'_v}^2)}{|f|^2\;
|\cosh(\theta)|^2\; \cosh^2(\theta)}\;,
\end{equation}
where $\theta$ is the function \eqref{theta_def}.

Thus the formulas for $K$ and $\varkappa$ related to the
representation \eqref{W2} are:
\begin{equation}\label{Kkappa_ftheta}
\begin{array}{lll}
K         &=& \Re \ds\frac{-({\theta'_u}^2+{\theta'_v}^2)}{|f|^2\; |\cosh(\theta)|^2\; \cosh^2(\theta)}\\[4ex]
\varkappa &=& \Im \ds\frac{-({\theta'_u}^2+{\theta'_v}^2)}{|f|^2\;
|\cosh(\theta)|^2\; \cosh^2(\theta)}\;.
\end{array}
\end{equation}
In order to express $K$ and $\varkappa$ by means of the functions
$w_j$, $(j=1;2)$ in the representation \eqref{W5}, we use
\eqref{thetapr_wpr} and \eqref{theta_w12}. Applying them to
\eqref{K+ikappa_ftheta} we find:
\begin{equation}\label{K+ikappa_fw12}
K+\ii\varkappa  = \ds\frac{- w'_1 \bar {w'_2}}{|f|^2\;
\left|\cosh\frac{w_1 + \bar w_2}{2}\right|^2\;
                  \cosh^2\frac{w_1 + \bar w_2}{2}}\;.
\end{equation}

The corresponding formulas in terms of $g_j$, $(j=1;2)$ in the
representation \eqref{W6} follow by using \eqref{thetapr_wpr},
\eqref{mf2_mf2} and \eqref{cosh^2theta_w12}:
\begin{equation*}%\label{K+ikappa_ftheta}
K+\ii\varkappa  = \ds\frac{- w'_1 \bar {w'_2}}{4|f|^2|g_1 g_2|\;
                  |\frac{1}{4}e^{-w_1}e^{-\bar w_2}(1+e^{w_1}e^{\bar w_2})^2|\;
                                     \frac{1}{4}e^{-w_1}e^{-\bar w_2}(1+e^{w_1}e^{\bar w_2})^2}\;.
\end{equation*}
By virtue of \eqref{g} и \eqref{wp}, the last formula takes the form:
\begin{equation}\label{K+ikappa_fg12}
K+\ii\varkappa  = \ds\frac{-4 g'_1 \bar {g'_2}}{|f|^2\;
                  |1 + g_1 \bar g_2|^2\;
                                    (1 + g_1 \bar g_2)^2}\;.
\end{equation}
Applying \eqref{E_fg1g2}, we get:
\begin{equation}\label{K+ikappa_Eg12}
K+\ii\varkappa  = \ds\frac{-4 g'_1 \bar {g'_2}}{E\; (1 + g_1 \bar
g_2)^2}\;.
\end{equation}

The corresponding formulas for $K$ and $\varkappa$, related to the
representation \eqref{W6} are:
\begin{equation}\label{Kkappa_fg1g2}
\begin{array}{lll}
K         &=& \Re \ds\frac{-4 g'_1 \bar {g'_2}}{|f|^2\;
                  |1 + g_1 \bar g_2|^2\; (1 + g_1 \bar g_2)^2}
                                    =\Re \ds\frac{-4 g'_1 \bar {g'_2}}{E\; (1 + g_1 \bar g_2)^2}\\[3ex]
\varkappa &=& \Im \ds\frac{-4 g'_1 \bar {g'_2}}{|f|^2\;
                  |1 + g_1 \bar g_2|^2\; (1 + g_1 \bar g_2)^2}
                                    =\Im \ds\frac{-4 g'_1 \bar {g'_2}}{E\; (1 + g_1 \bar g_2)^2}\;.
\end{array}
\end{equation}
The above formulas have been found by Asperti A. and Vilhena J. in
\cite{AA-VJ_DegGM}.

\section{The first fundamental form and the curvatures $K$, $\varkappa$, with respect to a canonical
Weierstrass representation.}
Let ${\M}$ be a minimal space-like surface of general type, parameterized by canonical
coordinates of the first type.

First we obtain a formula for the coefficient $E$ with respect to the canonical Weierstrass representation
\eqref{Wcanh}. Applying \eqref{can1w} and \eqref{thetapr_wpr}, we find the following formula for $|f|^2$:
\begin{equation}\label{mf^2_can_theta}
|f|^2 = \frac{1}{|w'_1 w'_2|} = \frac{1}{|{\theta'_u}^2+{\theta'_v}^2|}\;.
\end{equation}
Replacing into the general formula \eqref{E_ftheta}, we get:
\begin{equation}\label{E_Can_theta}
E = \ds\frac{|\cosh(\theta)|^2}{|{\theta'_u}^2+{\theta'_v}^2|}\;.
\end{equation}

To obtain a formula for $E$ with respect to the canonical Weierstrass representation \eqref{Wcanw},
we use \eqref{E_fw1w2} and \eqref{can1w}. Thus we find a formula for $E$ in terms of $w_1$ and $w_2$:
\begin{equation}\label{E_Can_w1w2}
E=\ds\frac{\left|\cosh\frac{w_1 + \bar w_2}{2}\right|^2}{|w'_1 w'_2|}.
\end{equation}

In a similar way, if ${\M}$ is given by \eqref{Wcang}, we replace $f$ into the general formula \eqref{E_fg1g2}
by the help of \eqref{can1g} and obtain:
\begin{equation}\label{E_Can_g1g2}
E=\ds\frac{|1+g_1\bar g_2|^2}{4|g'_1 g'_2|}.
\end{equation}

Next we find formulas for the curvatures $K$ and $\varkappa$. Using the representation \eqref{Wcanh},
we replace $f$ into the general formula \eqref{K+ikappa_ftheta} by means of \eqref{mf^2_can_theta},
and get:
\begin{equation}\label{K+ikappa_Can_theta}
K+\ii\varkappa  = \ds\frac{-|{\theta'_u}^2+{\theta'_v}^2|\, ({\theta'_u}^2+{\theta'_v}^2)}
                          {|\cosh(\theta)|^2\; \cosh^2(\theta)}\;.
\end{equation}

To obtain a formula for $K+\ii\varkappa$, when ${\M}$ is given by \eqref{Wcanw}, we replace into
the general formula \eqref{K+ikappa_fw12} the function $f$ by means of \eqref{can1w}.
Thus we have:
\begin{equation}\label{K+ikappa_Can_w12}
K+\ii\varkappa  = \ds\frac{- |w'_1 w'_2|\, w'_1 \bar {w'_2}}
{\left|\cosh\frac{w_1 + \bar w_2}{2}\right|^2\; \cosh^2\frac{w_1 + \bar w_2}{2}}\;.
\end{equation}
To obtain a formula for $K+\ii\varkappa$, when ${\M}$ is represented by \eqref{Wcang}, we replace into the general formula
\eqref{K+ikappa_fg12} the function $f$ by the help of \eqref{can1g}. Hence, we have:
\begin{equation}\label{K+ikappa_Can_g12}
K+\ii\varkappa  = \ds\frac{-16|g'_1 g'_2|\, g'_1 \bar {g'_2}}
                          {|1 + g_1 \bar g_2|^2\; (1 + g_1 \bar g_2)^2}\;.
\end{equation}
The formulas for the curvatures $K$ and $\varkappa$ with respect to the representation \eqref{Wcang} are
as follows:
\begin{equation}\label{Kkappa_Can_g12}
\begin{array}{lll}
K         &=& \Re \ds\frac{-16|g'_1 g'_2|\, g'_1 \bar {g'_2}}
                          {|1 + g_1 \bar g_2|^2\; (1 + g_1 \bar g_2)^2}\,,\\[4ex]
\varkappa &=& \Im \ds\frac{-16|g'_1 g'_2|\, g'_1 \bar {g'_2}}
                          {|1 + g_1 \bar g_2|^2\; (1 + g_1 \bar g_2)^2}\;.
\end{array}
\end{equation}

\section{Change of the functions $(g_1,g_2)$ under some basic geometric transformations of the minimal
space-like surface}

Let ${\M}$ be a minimal space-like surface of general type, parameterized by canonical coordi-nates $(u, v)$
of the first type. The complex variable $t$ is given by $t=u+\ii v$. We suppose that ${\M}$ is
given by the canonical representation \eqref{Wcang} by means of the pair $(g_1(t),g_2(t))$ of holomorphic functions.
The aim of this section is to study the changes of the pair $(g_1,g_2)$ under geometric transformations of
the surface .

First we consider the case of a motion of the surface ${\M}$ in $\RR^4_1$. We shall use some basic formulas
and facts about the spinors in $\RR^4_1$. Let us recall some of these formulas in a form useful for an application to
the theory of minimal space-like surfaces. To any vector $\x$ in $\RR^4_1$ we associate a Hermitian
$2\times 2$-matrix $S$ as follows:
\begin{equation}\label{Spin_S-x}
S=\left(
\begin{array}{rr}
     x_3 +x_4  &  \ii x_1 +x_2\\
-\ii x_1 +x_2  &     -x_3+x_4
\end{array}
\right) \ \leftrightarrow \
\x=(x_1,x_2,x_3,x_4)\;.
\end{equation}
This correspondence is a linear isomorphism between $\RR^4_1$ and the space of Hermitian $2\times 2$-matrices.
This correspondence has the following property: $\det S = -\x^2$. The last property means that from any linear
operator acting in the space of Hermitian $2\times 2$-matrices and preserving the determinant,
can be obtained an orthogonal operator in $\RR^4_1$.

If $\tilde A$ is a complex $2\times 2$-matrix, then $\tilde A S \tilde A^*$ is a Hermitian matrix,
where $\tilde A^*$ is the Hermitian conjugate of $\tilde A$. What is more, if $\det\tilde A = 1$, then
$\det \tilde A S \tilde A^* = \det S$. It follows from the above that to any matrix $\tilde A$ in
$\mathbf{SL}(2,\CC)$ corresponds an orthogonal matrix $A$ in $\mathbf{O}(3,1,\RR)$.
Therefore we have a group homomorphism $\tilde A \rightarrow A$, which can be written  as follows:
\begin{equation}\label{Spin_tildeA-A}
\hat S = \tilde A S \tilde A^* \ \rightarrow \ \hat\x = A\x\;.
\end{equation}
The so obtained homomorphism from $\mathbf{SL}(2,\CC)$ into $\mathbf{O}(3,1,\RR)$ is called \emph{spinor map}.
It is proved in the theory of spinors that the kernel of the spinor map consists of two elements: $\pm I$,
where $I$ is the unitary matrix. Further, the image of this map is the connected component of the unity element
in $\mathbf{O}(3,1,\RR)$, which is denoted in a standard way by $\mathbf{SO}^+(3,1,\RR)$.
Briefly speaking, this is the group of the matrices determining those transformations in $\RR^4_1$,
preserving not only the orientation of $\RR^4_1$, but also preserve both: the direction of time and
the orientation of the three-dimensional Euclidean subspace of $\RR^4_1$. These transformations of $\RR^4_1$
are called \emph{orthochronous} transformations.
The type of the kernel and the image of the spinor map \eqref{Spin_tildeA-A} implies that the spinor map induces
the following group isomorphism:
\begin{equation}\label{SL/+-1_SO+31}
\mathbf{SL}(2,\CC)/\{\pm I\}\ \cong \ \mathbf{SO}^+(3,1,\RR)\;.
\end{equation}
This means that $\mathbf{SL}(2,\CC)$ appears to be a two-sheeted covering of $\mathbf{SO}^+(3,1,\RR)$
and hence we can identify it with the spin group $\mathbf{Spin}(3,1)$ of $\mathbf{SO}^+(3,1,\RR)$.
In other words, \eqref{SL/+-1_SO+31} gives a representation of $\mathbf{Spin}(3,1)$ as $\mathbf{SL}(2,\CC)$,
which is called the spinor representation. Note that the group $\mathbf{SL}(2,\CC)$ is connected
and simply connected, and consequently it follows from the isomorphisms \eqref{Spin_tildeA-A}
and \eqref{SL/+-1_SO+31} that $\mathbf{SL}(2,\CC)$ also appears to be universal covering group
for $\mathbf{SO}^+(3,1,\RR)$.

Now, let $\x$ be an arbitrary complex vector in $\CC^4$. Up to now, considering different correspon\-dences,
we restricted $\x$ to be a real vector in $\RR^4_1$. It is an easy verification that the relations
\eqref{Spin_S-x} and \eqref{Spin_tildeA-A} are linear with respect to $\x$. Therefore, they are also valid
when $\x$ is an arbitrary complex vector in $\CC^4$. The only difference is that $S$ can be an arbitrary
(not necessarily Hermitian) complex matrix. Under a motion of the complex vector $\x$ with a matrix
in $\mathbf{SO}^+(3,1,\RR)$, the matrix $S$ is transformed in the same way, as it is described in
\eqref{Spin_tildeA-A}.

Let us return to minimal space-like surfaces. With the help of the above formulas we shall find how the
functions giving the Weierstrass representation of a minimal space-like surface are transformed
under a motion of the surface in $\RR^4_1$.

First, let the minimal space-like surface $({\M},\x)$ be parameterized by arbitrary isothermal coordinates.
If the surface $({\hat\M},\hat\x)$ is obtained from $({\M},\x)$ by means of orthochronous transformation
in $\RR^4_1$, then we have $\hat\x(t)=A\x(t)+\vb$, where $A \in \mathbf{SO}^+(3,1,\RR)$
and $\vb \in \RR^4_1$. The function $\Phi$ defined by \eqref{Phi_def}, as we noted by the formula
\eqref{hat_Phi-Phi-mov}, is transformed by: $\hat\Phi=A\Phi$. Next we introduce the complex matrix
$S_{\Phi}$, which is obtained by $\Phi$ according to the rule \eqref{Spin_S-x}:
\begin{equation}\label{SPhi_def}
S_{\Phi}=\left(
\begin{array}{rr}
     \phi_3 + \phi_4  &  \ii \phi_1 + \phi_2\\
-\ii \phi_1 + \phi_2  &     -\phi_3 + \phi_4
\end{array}
\right)\;.
\end{equation}
We denote by $\tilde A$ any of the two matrices in $\mathbf{SL}(2,\CC)$, corresponding to $A$ by means of the
homomorphism \eqref{Spin_tildeA-A}. If $S_{\hat\Phi}$ is the matrix obtained from $({\hat\M},\hat\x)$,
according to \eqref{Spin_tildeA-A} it is related to $S_{\Phi}$ as follows:
\begin{equation}\label{hatSPhi_SPhi}
S_{\hat\Phi} = \tilde A S_{\Phi} \tilde A^* \;.
\end{equation}

Now, suppose that $\M$ is given by a Weierstrass representation of the type \eqref{W6}.
By direct calculations we find:
\begin{equation*}%\label{phi_i1+2_3+4_fg12}
\begin{array}{l}
\phantom{-}\ii\phi_1+\phi_2=          -f(g_1 g_2+1)+f(g_1 g_2-1)=          -2f,\\
          -\ii\phi_1+\phi_2=\phantom{-}f(g_1 g_2+1)+f(g_1 g_2-1)=\phantom{-}2fg_1g_2,\\
\phantom{-}\phi_3+\phi_4=\phantom{-}f(g_1+g_2)+f(g_1-g_2)=\phantom{-}2fg_1,\\
          -\phi_3+\phi_4=          -f(g_1+g_2)+f(g_1-g_2)=          -2fg_2\;.
\end{array}
\end{equation*}
Consequently, the matrix $S_{\Phi}$ is represented by means of $f$, $g_1$ and $g_2$ as follows:
\begin{equation}\label{SPhi_fg12}
S_{\Phi}=\left(
\begin{array}{ll}
     2fg_1  &  -2f \\
  2fg_1g_2  &  -2fg_2
\end{array}
\right)\;.
\end{equation}
Denoting the elements of $S_{\Phi}$ by $s_{ij}$, then we have the following expressions for  $f$, $g_1$ and $g_2$:
\begin{equation}\label{fg12_s1234}
f=-\frac{1}{2}s_{12}, \quad g_1=-\frac{s_{11}}{s_{12}}, \quad g_2=\frac{s_{22}}{s_{12}}\;.
\end{equation}

Since $S_{\Phi}$ is transformed according the rule \eqref{hatSPhi_SPhi}, then via \eqref{fg12_s1234}
we shall find the transfor-mation formulas for the functions $f$, $g_1$ and $g_2$. For that purpose we
denote the elements of $\tilde A$ in the following way:
\begin{equation}\label{tildeA_abcd}
\tilde A=\left(
\begin{array}{rr}
     \bar a  & -\bar b \\
    -\bar c  &  \bar d
\end{array}
\right)\,;
\quad a,b,c,d\in \CC\,;
\quad ad-bc=1\;.
\end{equation}
After multiplying the matrices in \eqref{hatSPhi_SPhi} and simplifying, we get:
\begin{equation}\label{hatSPhi_fg12}
S_{\hat\Phi}=\left(
\begin{array}{ll}
     2f(a g_1+b)(         - \bar b g_2 + \bar a)   &  2f(c g_1+d)(\phantom{-}\bar b g_2 - \bar a) \\
     2f(a g_1+b)(\phantom{-}\bar d g_2 - \bar c)   &  2f(c g_1+d)(         - \bar d g_2 + \bar c)
\end{array}
\right)\;.
\end{equation}
Applying \eqref{fg12_s1234} to $\hat f$, $\hat g_1$ and $\hat g_2$, we find the transformation formulas
of the functions in the Weierstrass representation of the type \eqref{W6} under an orthochronous transformation
of $\M$ in $\RR^4_1$:
\begin{equation}\label{hatfg12_fg12}
\begin{array}{l}
\hat f = f(cg_1+d)(-\bar b g_2 + \bar a)\,;\\[0.7ex]
\hat g_1 = \ds\frac{ag_1+b}{cg_1+d}\,; \quad \hat g_2 = \ds\frac{\phantom{-}\bar d g_2 - \bar c}{-\bar b g_2 + \bar a}\;.
\end{array}
\end{equation}

Now, let us consider the inverse statement. Suppose that $({\hat\M},\hat\x)$ and $({\M},\x)$
are two minimal space-like surfaces, given by the Weierstrass representation of the type \eqref{W6}
related by means of \eqref{hatfg12_fg12}. We shall show that they can be obtained one from the other by
an orthochronous transformation in $\RR^4_1$. For that purpose, we introduce $\tilde A$ by means of \eqref{tildeA_abcd}.
Let $A$ be the corresponding to $\tilde A$ matrix under the homomorphism \eqref{Spin_tildeA-A}.
With the help of $A$ we obtain a third surface $({\hat{\hat\M}},\hat{\hat\x})$ given by the formula: $\hat{\hat\x} = A\x$.
We proved that $\hat{\hat\M}$ has a Weierstrass representation with functions also satisfying \eqref{hatfg12_fg12}.
Therefore, $\hat\M$ and $\hat{\hat\M}$ are generated by one and the same functions by means of formulas \eqref{W6}
and consequently they are obtained one from the other by a translation in $\RR^4_1$. Since $\hat{\hat\M}$ is obtained from
$\M$ by an orthochronous transformation, then $\hat\M$ is also obtained from $\M$ by orthochronous transformation.
Summarizing we obtain the following statement:
\begin{thm}\label{A-B_W6}
Let $({\hat\M},\hat\x)$ and $({\M},\x)$ be two minimal space-like surfaces in $\RR^4_1$, given by Weierstrass
representations of the type \eqref{W6}. The following conditions are equivalent:
\begin{enumerate}
    \item $({\hat\M},\hat\x)$ and $({\M},\x)$ are related by an orthochronous transformation in $\RR^4_1$ of the type:\\
    $\hat\x(t)=A\x(t)+\vb$, where $A \in \mathbf{SO}^+(3,1,\RR)$ and $\vb \in \RR^4_1$.
    \item The functions in the Weierstrass representations of $({\hat\M},\hat\x)$ and $({\M},\x)$ are related by equalities
    \eqref{hatfg12_fg12}, where $a,b,c,d\in \CC$, $ad-bc=1$.
\end{enumerate}
\end{thm}
 Up to now we considered only the case of a motion from the connected component of the identity in $\mathbf{O}(3,1,\RR)$.
Next we show that any of the three remaining cases can be reduced to the considered one. Let us consider the case
of a transformation, which is not orthochronous. Such a concrete transformation can be obtained by a change
of the signs of the four coordinates: $\hat\x(t)=-\x(t)$. This implies the change of the sign of the function $f$
in the Weierstrass representation \eqref{W6}, while the functions $g_1$ and $g_2$ remain the same.
Any non-orthochronous transformation can be obtained as a composition of this concrete transformation and an orthochronous
transformation in $\RR^4_1$. Therefore, if two minimal space-like surfaces are obtained one from the other by a
non-orthochronous transformation, then the functions in the Weierstrass representation are related by formulas,
which are similar to \eqref{hatfg12_fg12} with the only difference in the sign of the formula for $\hat f$.
Further, we consider the case of a non-orthochronous improper transformation. An example of such a transformation is
the symmetry with respect to the hyperplane $x_4=0$ which is given by the change of the sign of $x_4$.
This implies a change of the places of both functions $g_1$ and $g_2$ in the Weierstrass representation, while the function
$f$ remains the same. Any non-orthochronous improper transformation can be obtained as a composition of this symmetry and
an orthochronous transformation. Therefore the functions in the Weierstrass representation are changed similarly to
\eqref{hatfg12_fg12}, but this time the formulas for $\hat g_1$ and $\hat g_2$ change their places, while the formula
for $\hat f$ is the same. Finally, we consider the case of an orthochronous improper transformation. Such a transformation
can be obtained as a combination of the last two cases. Therefore, the transformation formulas for the functions
in the Weierstrass representation are obtained from \eqref{hatfg12_fg12}, by the change of the sign of $\hat f$
and the change of the places of the formulas for $\hat g_1$ and $\hat g_2$.

Now, let $\hat\M$ and $\M$ be two minimal space-like surfaces, parameterized by canonical coordinates
and the surface $\hat\M$ is obtained from $\M$ by a motion in $\RR^4_1$. Suppose that $\M$
is given by a canonical Weierstrass representation of the type \eqref{Wcang}. Since $\hat\x=A\x+\vb$
implies that $\hat\Phi'=A\Phi'$, then we have $\hat\Phi'{}^2=\Phi'^2=1$. Consequently the canonical coordinates
of $\M$ appear to be also canonical coordinates for $\hat\M$. Taking into account that the canonical
Weierstrass representation \eqref{Wcang} is a special case of the representation \eqref{W6}, then
the pair $(g_1,g_2)$ is transformed by the formulas \eqref{hatfg12_fg12}. Note that these formulas can be applied
in the cases of an orthochronous or a non-orthochronous transformation. This is so because the two cases
differ from each other only the formula for the function $f$. Further, we see that it only remain the linear
fractional functions from \eqref{hatfg12_fg12}, which allowas us to replace the condition $ad-bc=1$ with the
more general condition $ad-bc\neq 0$. This is possible, because the linear fractional function does not change
if its matrix is multiplied by a non zero factor.

Summarizing the above remarks, in view of Theorem \ref{A-B_W6} we obtain the following statement.
\begin{thm}\label{A-B_Cang}
Let $({\hat\M},\hat\x)$ and $({\M},\x)$ be two minimal space-like surfaces of general type, given by
the canonical Weierstrass representation of the type \eqref{Wcang}. The following conditions are equivalent:
\begin{enumerate}
    \item $({\hat\M},\hat\x)$ and $({\M},\x)$ are related by a transformation in $\RR^4_1$ of the type:\\
    $\hat\x(t)=A\x(t)+\vb$, where $A \in \mathbf{SO}(3,1,\RR)$ and $\vb \in \RR^4_1$.
    \item The functions in the Weierstrass representations of $({\hat\M},\hat\x)$ and $({\M},\x)$ are related by the
    following equalities:
\begin{equation}\label{hatg12_g12}
\hat g_1 = \ds\frac{ag_1+b}{cg_1+d}\,; \quad \hat g_2 = \ds\frac{\phantom{-}\bar d g_2 - \bar c}{-\bar b g_2 + \bar a}\;,
\end{equation}
    where $a,b,c,d\in \CC$, $ad-bc\neq 0$.
\end{enumerate}
\end{thm}

If $({\hat\M},\hat\x)$ and $({\M},\x)$ are related by an improper transformation in $\RR^4_1$, then
in the formulas \eqref{hatg12_g12} one has to change the places of $\hat g_1$ and $\hat g_2$.

In the end, we write down \eqref{hatg12_g12} in a form, which is useful for applications.
For that purpose, let us denote by $Gz$, where $G\in \mathbf{GL}(2,\CC)$ and $z\in \CC$, the standard action
of the group $\mathbf{GL}(2,\CC)$ in the complex plane by means of linear fractional transformations.
Denoting by $B$ the matrix of the linear fractional function for $\hat g_1$ in \eqref{hatg12_g12},
by direct computations we see that the matrix of $\hat g_2$ is up to a factor the matrix ${B^*}^{-1}$.
Hence, the formulas \eqref{hatg12_g12} can be written briefly as follows:
\begin{equation}\label{matr_hatg12_g12}
\hat g_1 = Bg_1\,; \quad \hat g_2 = {B^*}^{-1} g_2\,;   \quad B\in \mathbf{GL}(2,\CC)\;.
\end{equation}

Finally we give a natural approach to the family of the minimal space-like surfaces of general type,
associated with a given one. Let $(g_1(t),g_2(t))$ be a pair of holomorphic functions defined in
a disc $\mathcal D$, centered at $(0, 0)$ in the parametric plane $\CC$. Consider the minimal
space-like surface $({\M},\x)$, generated by the pair $(g_1(t),g_2(t))$ by means of \eqref{Wcang}.
For any complex number $a, \, |a|=1$ we introduce the pair of holomorphic functions
\begin{equation}\label{g1g2_as+b}
(\tilde{g}_1(s),\tilde{g}_2(s))=(g_1(as),g_2(as)); \quad s \in \mathcal D
\end{equation}
and denote by ${\tilde{\M}}$ the minimal space-like surface, generated by the pair
$(\tilde{g}_1,\tilde{g}_2)$ by means of \eqref{Wcang}. Further we denote by $\Phi$, $\Psi$ and
$\tilde\Phi$, $\tilde\Psi$ the corresponding vector holomorphic functions on ${\M}$ and ${\tilde{\M}}$.
Replacing \eqref{g1g2_as+b} into the representation \eqref{Wcang}, we get the following relation
between $\tilde\Phi$ and $\Phi$:
\begin{equation}\label{Phi_as+b}
\tilde\Phi(s)=\ds\frac{1}{a}\:\Phi(as).
\end{equation}
Since $\tilde \Phi'^2=1$, then $s=\frac{t}{a}$ determines canonical coordinates on $\tilde{\M}$.

After an integration we obtain the corresponding formula for $\tilde\Psi$:
\begin{equation}\label{Psi_as+b}
\tilde\Psi(s)=\ds\frac{1}{a^2}\:\Psi(as).
\end{equation}

Denoting $a=e^{\ii\frac{\varphi}{2}}$ and $\tilde{\M}=\M_{\varphi}$, we have:

Any minimal space-like surface of general type $\M : \x=\x(t); \; t\in \mathcal D$ generates a one-parameter
family $\{\M_{\varphi}\}$ of minimal space-like surfaces, given by the formula
$$\M_\varphi: \;  \x_\varphi (s) = \Re(e^{-\ii\varphi}\:\Psi(e^{\ii \frac{\varphi}{2}}\, s)); \quad \varphi \in
[0, \frac{\pi}{2}], \quad s \in \mathcal D, $$
where $s=e^{-\ii \frac{\varphi}{2}}t$ determines canonical coordinates on $\M_\varphi$.

The surfaces of the family $\{\M_\varphi\}$ are said to be \emph{associated} with the given surface $\M$.

Since the generating  holomorphic functions of the family of the associated surfaces are given by
\eqref{g1g2_as+b}, taking into account formulas \eqref{E_Can_g1g2} and \eqref{Kkappa_Can_g12}, we observe that
the transformation $\M \, \rightarrow \,\M_{\varphi}$ given by $s \, \rightarrow \, e^{\ii\frac{\varphi}{2}}s$
preserves $E$, $K$ and $\varkappa$, i.e. it is a \emph{special isometry between $\M$ and $\M_{\varphi}$ preserving
the normal curvature $\varkappa$}.

 Denote by ${\bar \M}$ the minimal space-like surface conjugate to ${\M}$, which is given by the formula
$\y=\Im (\Psi) = \Re (-\ii\Psi)$. Then ${\bar \M}$ is the associated with $\M$ surface $\M_{\varphi},\;
\varphi = \frac{\pi}{2}$.

Thus we have:

If the minimal space-like surface $\M$, parameterized by canonical coordinates, is generated by the pair
$(g_1(t),g_2(t))$, then the minimal space-like surface ${\bar \M}$, conjugate to $\M$, is generated by the pair
$(g_1(e^{\ii \frac{\pi}{4}}\, s),g_2(e^{\ii \frac{\pi}{4}}\, s))$ with canonical parameter
$s=e^{-\ii \frac{\pi}{4}}\, t$.

If $t$ determines canonical coordinates of the first type on ${\M}$, then $e^{-\ii \frac{\pi}{4}}\, t$
gives canonical coordinates of the second type on $\M$ and vice versa.

Hence:

\emph{The canonical coordinates of the second type on $\M$ are canonical coordinates of the first type on $\bar{\M}$
and the canonical coordinates of the first type on ${\M}$ are canonical coordinates of the second type
on ${\bar \M}$.}
\vskip 3mm
\textbf{Acknowledgments:}
The first author is  partially supported by the National Science Fund,
Ministry of Education and Science of Bulgaria under contract DFNI-I 02/14.


\begin{thebibliography}{99}


\bibitem{A-P-1}
Al\'{\i}as L., Palmer, B., \emph{Curvature properties of zero mean curvature surfaces in four-dimensional
{Lorenzian} space forms.}, Mathematical Proceedings of the Cambridge Philosophical Society,
\textbf{124} (1998), 315-327.

\bibitem{AA-VJ_DegGM}
Asperti A., Vilhena J., \emph{Spacelike Surfaces in $\mathbb{L}^4$ with Degenerate Gauss Map.},
Results in Mathematics, \textbf{60} (2011), 185-211.

\bibitem{AA-VJ_PrescGM}
Asperti A., Vilhena J., \emph{Spacelike Surfaces in $\mathbb{L}^4$ with prescribed Gauss
map and nonzero mean curvature.}, Matem\'{a}tica Contempor\^{a}nea, \textbf{33} (2007), 55–83.

\bibitem{E-1}
Eisenhart L., \emph{A Fundamental Parametric Representation of Space Curves},
Ann. Math., Second Series, \textbf{13} (1/4) (1911 - 1912), 17-35.

\bibitem{E-R-1}
Estudillo F. J. M., Romero A., \emph{On maximal surfaces in the n-dimensional Lorentz-Minkowski space.},
Geom. Dedicata, \textbf{38} (1991), 167-174.

\bibitem{G-2}
Ganchev G., \emph{Canonical Weierstrass Representation of Minimal and Maximal Surfaces
in the Three-dimensional Minkowski Space.}
arXiv:0802.2632 (https://arxiv.org/abs/0802.2632)

\bibitem{G-K-2}
Ganchev G., Kanchev K., \emph{Canonical Weierstrass representations for minimal surfaces in Euclidean 4-space},
arXiv:1609.01606 (https://arxiv.org/abs/1609.01606)

\bibitem{HO}
Hoffman D., Osserman R., \emph{The geometry of the generalized Gauss map},
Memoirs of the American Mathematical Society, \textbf{28} (236) (1980).

\bibitem{dM}
Montcheuil M. de, \emph{R\'{e}solution de l\textquoteright\'{e}quation $ds^2=dx^2+dy^2+dz^2$.},
Bulletin de la Soci\'{e}t\'{e} Math\'{e}matique de France, \textbf{33} (1905), 170-171.

\end{thebibliography}
\end{document}